\documentclass[final,3p,times,onecolumn]{elsarticle}
\usepackage{subfig}
\usepackage{amsmath}
\usepackage{upgreek}
\usepackage{multirow}
\usepackage{siunitx}
\usepackage[export]{adjustbox}
\usepackage[linesnumbered,ruled, noend]{algorithm2e}

\SetCommentSty{mycommfont}

\usepackage{booktabs}
\usepackage{scalerel}

\usepackage[hidelinks]{hyperref}
\usepackage[nameinlink, capitalise]{cleveref}

\Crefname{figure}{Fig.}{Figs.}
\Crefname{equation}{Eq.}{Eqs.}
\makeatletter
\def\ps@pprintTitle{%
  \let\@oddhead\@empty
  \let\@evenhead\@empty
  \let\@oddfoot\@empty
  \let\@evenfoot\@oddfoot
}
\makeatother

\begin{document}

\begin{frontmatter}


\title{Surrogate-based multiscale analysis of experiments on thermoplastic composites under off-axis loading}
\author[a]{M. A. Maia\corref{cor1}}
\author[a]{I. B. C. M. Rocha}
\author[a]{D. Kovacevic}
\author[a]{F. P. van der Meer}

\address[a]{{Delft University of Technology, Department of Civil Engineering and Geosciences}, {P.O.Box 5048}, 2600GA, Delft, {The Netherlands}}
\cortext[cor1]{Corresponding author. E-mail: marinaalvesmaia@gmail.com}

\begin{abstract}
In this paper, we present a surrogate-based multiscale approach to model constant strain-rate and creep experiments on unidirectional thermoplastic composites under off-axis loading. In previous contributions, these experiments were modeled through a single-scale micromechanical simulation under the assumption of macroscopic homogeneity. Although efficient and accurate in many scenarios, simulations with low-off axis angles showed significant discrepancies with the experiments. It was hypothesized that the mismatch was caused by macroscopic inhomogeneity, which would require a multiscale approach to capture it. However, full-field multiscale simulations remain computationally prohibitive. To address this issue, we replace the micromodel with a Physically Recurrent Neural Network (PRNN), a surrogate model that combines data-driven components with embedded constitutive models to capture history-dependent behavior naturally. The explainability of the latent space of this network is also explored in a transfer learning strategy that requires no re-training. With the surrogate-based simulations, we confirm the hypothesis raised on the inhomogeneity of the macroscopic strain field and gain insights into the influence of adjustment of the experimental setup with oblique end-tabs. Results from the surrogate-based multiscale approach show better agreement with experiments than the single-scale micromechanical approach over a wide range of settings, although with limited accuracy on the creep experiments, where macroscopic test effects were implicitly taken into account in the material properties calibration. 
\end{abstract}


\end{frontmatter}


\newcommand{\fesqr}{FE$^2$}
\newcommand{\pd}[2]{\frac{\partial{#1}}{\partial{#2}}}
\newcommand{\strainvec}{\boldsymbol{\varepsilon}}
\newcommand{\stressrot}{\boldsymbol{\sigma}_{\textrm{F}}}
\newcommand{\predstresslocal}{\widehat{\boldsymbol{\sigma}}^{\Omega}_{\textrm{L}}}
\newcommand{\predstressrot}{\widehat{\boldsymbol{\sigma}}^{\Omega}_{\textrm{F}}}
\newcommand{\predstress}{\widehat{\boldsymbol{\sigma}}^{\Omega}}
\newcommand{\strainstresscurve}{\boldsymbol{\varepsilon}\text{-}\boldsymbol{\sigma}}
\newcommand{\intvarvec}{\boldsymbol{\alpha}}
\newcommand{\macrostrainvec}{\boldsymbol{\varepsilon}^{\Omega}}
\newcommand{\macrostressvec}{\boldsymbol{\sigma}^{\Omega}}
\newcommand{\microstressvec}{\boldsymbol{\sigma}^{\omega}}
\newcommand{\microstrainvec}{\boldsymbol{\varepsilon}^{\omega}}
\newcommand{\macrodispfield}{\mathbf{u}^\Omega}
\newcommand{\microdispfield}{\mathbf{u}^\omega}
\newcommand{\dispfield}{\mathbf{u}}
\newcommand{\weights}{\mathbf{W}}
\newcommand{\Wjmat}{\mathbf{W}_{j}}
\newcommand{\Ljmat}{\mathbf{W}_j}
\newcommand{\biases}{\mathbf{b}}
\newcommand{\hidden}{\mathbf{h}}
\newcommand{\histunits}{\mathbf{h}}
\newcommand{\stiffmatrix}{\mathbf{D}}
\newcommand{\xvec}{\mathbf{x}}
\newcommand{\Bmat}{\mathbf{B}}
\newcommand{\macrodmat}{\mathbf{D}^{\Omega}}
\newcommand{\Ftensor}{\mathbf{F}}
\newcommand{\Rtensor}{\mathbf{R}}
\newcommand{\Utensor}{\mathbf{U}}
\newcommand{\predstressunrot}{\widehat{\boldsymbol{\sigma}}_\text{U}^{\Omega}}
\newcommand{\stressunrot}{\boldsymbol{\sigma}_\text{U}}
\newcommand{\Umacrotensor}{\mathbf{U}^{\Omega}}
\newcommand{\Umicrotensor}{\mathbf{U}^{\omega}}
\newcommand{\Bmacrotensor}{\mathbf{B}^{\Omega}}
\newcommand{\Bmicrotensor}{\mathbf{B}^{\omega}}
\newcommand{\Idtensor}{\mathbf{I}}
\newcommand{\Fmacrotensor}{\mathbf{F}^{\Omega}}
\newcommand{\Fmicrotensor}{\mathbf{F}^{\omega}}
\newcommand{\Rmacrotensor}{\mathbf{R}^{\Omega}}
\newcommand{\Rmacrolocal}{\mathbf{R}^{\Omega}_{\textrm{L}}}
\newcommand{\Umacrolocal}{\mathbf{U}^{\Omega}_{\textrm{L}}}
\newcommand{\Fmacrolocal}{\mathbf{F}^{\Omega}_{\textrm{L}}}
\newcommand{\Qmat}{\mathbf{Q}_{0}}
\newcommand{\incdefgrad}{\Delta\mathbf{F}}
\newcommand{\defgrad}[2][]{\mathbf{F}_{#1}^{#2}}
\newcommand{\stress}[2][]{\boldsymbol{\sigma}_{#1}^{#2}}
\newcommand{\intvar}[2][]{\boldsymbol{\alpha}_{#1}^{#2}}
\newcommand{\predmacrostress}{\widehat{\boldsymbol{\sigma}}^{\Omega}_t}
\newcommand{\norm}[1]{\left\lVert#1\right\rVert}

% Definition of datasets
\newcommand{\testmono}{\mathcal{T}_\text{mono}^{\textrm{16\,modes}}}
\newcommand{\testslower}{\mathcal{T}_{\text{slower}}}
\newcommand{\testfaster}{\mathcal{T}_{\text{faster}}}
\newcommand{\testunlfixed}{\mathcal{T}_{\text{unl}}^{\text{fixed}}}
\newcommand{\testunlgpbased}{\mathcal{T}_{\text{16}}}
\newcommand{\testunlgp}{\mathcal{T}_{\text{unl}}^{\text{non--prop.\,GP}}}

% Definition of properties sets
\newcommand{\material}[1]{$\mathcal{C}_{\textrm{#1}}^{\omega}$}
\newcommand{\setcsr}[1]{$\boldsymbol{\theta}_{\textrm{#1}}^{\dot{\varepsilon}_\textrm{yy}^{\textrm{eng}}}$}
\newcommand{\setcreep}[1]{$\boldsymbol{\theta}_{\textrm{#1}}^{\textrm{creep}}$}
\newcommand{\setcsrmode}[2]{$\boldsymbol{\theta}_{\textrm{#1}}^{\dot{\varepsilon}_\textrm{yy}^{\textrm{eng}},\,#2 \,\textrm{mode}}$}
\newcommand{\setcsrmodes}[2]{$\boldsymbol{\theta}_{\textrm{#1}}^{\dot{\varepsilon}_\textrm{yy}^{\textrm{eng}},\,#2 \,\textrm{modes}}$}
\newcommand{\setcreepmode}[2]{$\boldsymbol{\theta}_{\textrm{#1}}^{\textrm{creep},\,#2 \,\textrm{mode}}$}
\newcommand{\setcreepmodes}[2]{$\boldsymbol{\theta}_{\textrm{#1}}^{\textrm{creep},\,#2 \,\textrm{modes}}$}

\section{Introduction}
\label{sec:intro}

In many engineering applications, understanding how phenomena at lower scales affect the material behavior at the macroscale is essential. For composite materials, this knowledge can lead to improved designs at lower costs. To achieve this, systematic experiments are needed to characterize the material and the response of the designed structural component, processes that can be costly and time-consuming. Fortunately, with advances in computational power, numerical models based on finite element (FE) can help alleviate the costs of the experimental campaigns by replacing part of these tests with computational simulations. A useful computational model, however, still relies on representative constitutive models and material properties. 

In that regard, several methods have been proposed to predict the constitutive behavior of complex heterogeneous materials, including composites, with varying degrees of accuracy, generality, and efficiency \cite{MATOUS2017}. Some can be highly efficient but are also bound to the assumptions made in their formulation, restricting their flexibility. For example, mean-field homogenization schemes (e.g. Mori-Tanaka \cite{MORI1973} and Self-Consistent method) based on Eshelby's solution \cite{Eshelby1957} are rather limited regarding the geometric features of the inclusions, volume fraction, stress localization and material nonlinearity. We refer to \cite{JAIN2019,HESSMAN2021} for recent reviews on this type of solution. Another efficient way to predict the homogenized constitutive response of the material is to use mesoscale (ply-level) models. These models have a certain number of parameters that should be calibrated based on experimental data or numerical micromodel simulations. While fast, the balance between the number of parameters, data availability, and ability to predict general loading conditions is not trivial. Mesoscale models can be seriously limited if certain assumptions and interactions are not considered. This is illustrated in \cite{VANDERMEER2016}, where a study on a recently developed homogenized orthotropic plasticity model quantified the loss of accuracy due to necessary simplifications such as ignoring the influence of stress in the fiber direction on the plasticity and use of a constant plastic Poisson ratio.

At the opposite end of the spectrum in terms of efficiency are multiscale methods, particularly FE$^2$. This approach associates the constitutive response of an integration point at the macroscale with the homogenized response of the micromodel, or Representative Volume Element (RVE), leveraging the convenience of using simpler constitutive models at the lower scale of observation. At this level, the explicit microstructure description offers more information on material integrity, and phenomena such as orthotropic behavior, (visco-)plasticity, aging, damage and strain/stress localization can be more easily incorporated. This generality, however, is associated with extreme computational cost, limiting its feasibility in practical applications. 

Numerous approaches to improve the computational efficiency of FE$^2$ have been proposed (e.g. parallel computing, fast Fourier transform methods \cite{DEGEUS2017}, clustering-based methods \cite{FERREIRA2022,Deng2022,CHAOUCH2024}, etc.), but here we focus on the ones based on machine learning techniques. Specifically, frameworks built on top of surrogate models that replace the micromodel problem entirely. A surrogate model typically approximates the mapping between the homogenized quantities and usually consists of a data-driven model trained on snapshots of the micromodel, or in the case of history-dependent materials, sequences of snapshots. Examples include Gaussian Processes (GPs) \cite{ROCHA2021,FUHG2022}, Deep Neural Networks \cite{Eivazi2023,FENG2022,Aldakheel2023}, Graph Neural Networks \cite{VLASSIS2023}, Transformers \cite{Zhongbo2024,Pitz2024}, Convolution Neural Networks \cite{ABUEIDDA2021}, Recurrent Neural Networks (RNNs) \cite{WU2022114476, LOGARZO2021, GORJI2020, Mozaffar26414,ghane2023recurrent,CHEUNG2024}, frameworks that combine some of these methods with Reduced Order Modeling techniques \cite{GUO2018,GUO2019,CASENAVE2020} and, more recently, hybrid models, which combine machine learning techniques with some knowledge of the problem to fulfill thermodynamics and kinematic constraints by construction. Despite the wide popularity in the literature, the purely data-driven models struggle in the extrapolation regime and are often data-hungry, and hybrid models remain specialized to a narrow range of constitutive behavior (e.g. homogeneous and hyperelastic materials \cite{Klein2023,LINKA2023}).

To address those issues, we proposed the Physically Recurrent Neural Networks (PRNNs) \cite{Maiaetal2023}. The idea is to preserve the knowledge we have on constitutive models, specifically the ones employed in the RVE used to generate the high-fidelity data and train the network, by embedding them in an encoder-decoder architecture. This combination results in great generalization properties with small training sets and offers a clear interpretation of its latent space. In our latest development \cite{Maiaetal2024}, we extended the applicability of PRNNs to deal with rate and path-dependent 3D problems in a finite strain formulation. The inclusion of damage models in the PRNN has also been recently addressed in \cite{kovacs2024} for the microscale debonding of composite materials. 

Our starting point are the experimental results and contributions in \cite{Kovacevic2022,Kovacevic2022exp,Kovacevic2024creep}. In these works, constant strain-rate and creep experiments were modeled at the microscale under the assumption that macroscale variations are negligible. In \cite{Kovacevic2022}, the authors laid the foundations of a custom framework based on a single-scale micromodel simulation with special boundary conditions to emulate a uniaxial stress state in the global frame for a composite with rate-dependent response under arbitrary off-axis loading. Later, the formulation was expanded to include failure analysis for constant strain-rate \cite{Kovacevic2022exp}, constant stress \cite{Kovacevic2024creep} and cyclic loading \cite{Kovacevic2024uni}. For different strain-rates, stress levels and off-axis angles, the micromodel-based formulation was able to accurately reproduce the experimental results in most cases, with the exception of small off-axis angle scenarios, in which a significant difference was observed already prior to failure. This difference suggests a violation of the assumption of macroscopic uniform deformations. However, testing this hypothesis with the same high-fidelity micromodel would require a multiscale approach, which is computationally intractable. Now with the PRNN developed in \cite{Maiaetal2024}, a surrogate-based multiscale analysis comes within reach. 

Other works that explored the surrogate-based multiscale approach to predict the off-axis homogenized response of composites are \cite{Gao2020} and \cite{Kalina2023}. \citet{Gao2020} validated the results against experimental data but did not extend the formulation to account for different strain-rates and temperatures, while \citet{Kalina2023} focused on the numerical modeling aspects considering hyperelastic models only. Lastly, the strain inhomogeneity matter was also discussed in \citet{Yifan2023} through the lens of the Digital Image Correlation method. In that case, the authors focused on the characterization of the material for different experimental setups, leaving out details on the computational modeling. 

In this contribution, we model the constant strain-rate and creep experiments in \cite {Kovacevic2022exp,Kovacevic2024creep} as a multiscale problem with the PRNN as a surrogate model to the homogenized response of a micromodel with history-dependent materials. With this, we contribute to a better understanding of the experiments by testing the hypotheses raised in \cite{Kovacevic2022exp} and \cite{Kovacevic2024creep} on the macroscopic strain/stress variability. In addition to that, we illustrate new features of the PRNN for transfer learning, favoring more efficient training and online evaluation. In the following, we discuss the test setup and methods. Then, we assess the performance of the surrogate-based approach on a selection of experiments, showcasing the robustness of the network for multiscale analysis. Finally, we draw the main conclusions of this work in \Cref{sec:conclusion}.

\section{Methods}
\label{sec:fullmulti}
\begin{figure}[!ht]
\centering
\hfill
\subfloat[Geometry with dimensions in mm\label{fig:coupongeometry}]{\includegraphics[width=.255\textwidth]{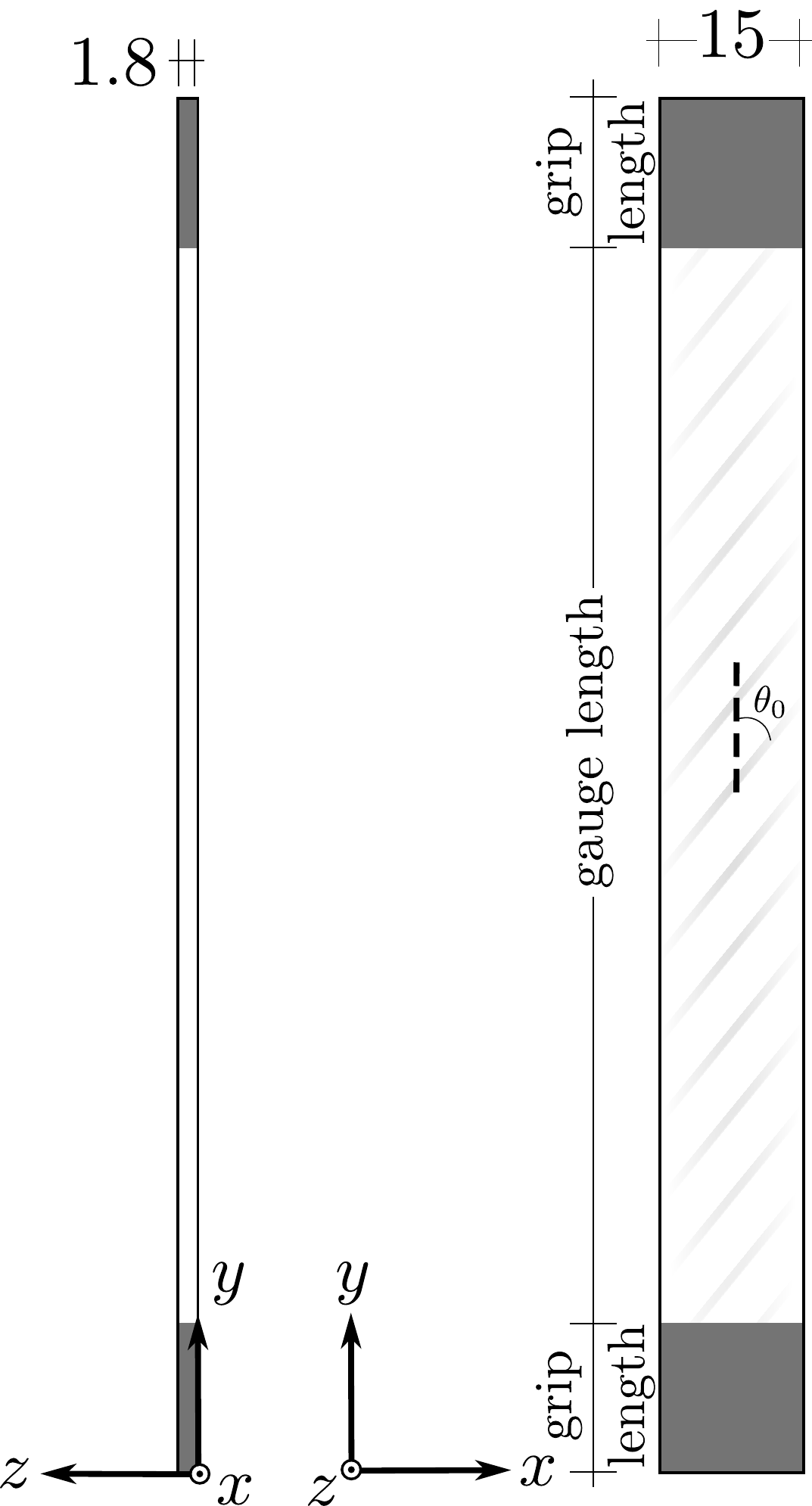}}
\hfill
\subfloat[Setup used for creep testing \cite{Kovacevic2024creep}\label{fig:creepmachine}]{\includegraphics[width=.255\textwidth]{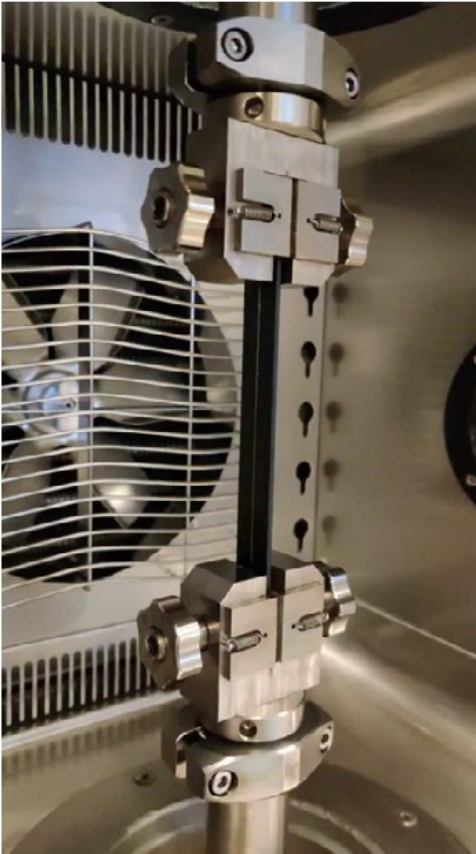}}
\hfill \
\caption{Schematic representation of UD carbon/PEEK composite system specimen and creep testing setup.}
\label{fig:specimen}
\end{figure}

In this section, we outline the most relevant aspects of the experimental and numerical setups used to obtain the strain-stress response of a thermoplastic composite system over different off-axis angles, strain-rates and stress levels. 

\subsection{Experimental setup}
\label{subsec:experimental}

Coupons were made of carbon fiber reinforced unidirectional tapes and PEEK matrix with fiber volume fraction of 0.4. Two types of experiments are considered: constant strain-rate and creep. The geometry of the coupons manufactured for both experiments is illustrated in \Cref{fig:coupongeometry}, where $\theta_0$ is the off-axis angle. The numerical counterpart models only the gauge length, which is \qty{120}{\mm} for the constant strain rate experiments. In these tests, the crosshead speeds were kept constant corresponding to strain-rates ranging between $10^{-6}$ \unit{\per\s} and $10^{-3}$ \unit{\per\s}. To record the stress-strain relationship, clip-on expensometers were used. Crosshead displacements were converted to engineering strain, while the resulting force was used to compute the engineering stress. 

For the creep experiments, the gauge length was \qty{100}{\mm} for the off-axis angles \qty{90}{\degree} and \qty{45}{\degree}, and \qty{120}{\mm} for the smaller angles \qty{30}{\degree} and \qty{15}{\degree}. In these experiments, constant force was applied on the coupon, from which the engineering stress was calculated and the crosshead displacements gave the engineering strain. Further details on the experimental setup can be found in \cite{Kovacevic2022exp,Kovacevic2024creep}. 

\subsection{Multiscale problem formulation}
\label{subsec:fe2}

To replicate the experiments discussed in the previous section, a concurrent multiscale approach (FE$^2$) can be used. In this case, the macroscopic domain $\Omega$ is discretized into a Finite Element (FE) mesh, with a periodic representative volume element (RVE) $\omega$ nested to each integration point. The RVE consists of another FE model that characterises the heterogeneous material at a length scale significantly lower than the macroscopic one.
\begin{figure}[!h]
\centering
\includegraphics[width=.6\textwidth]{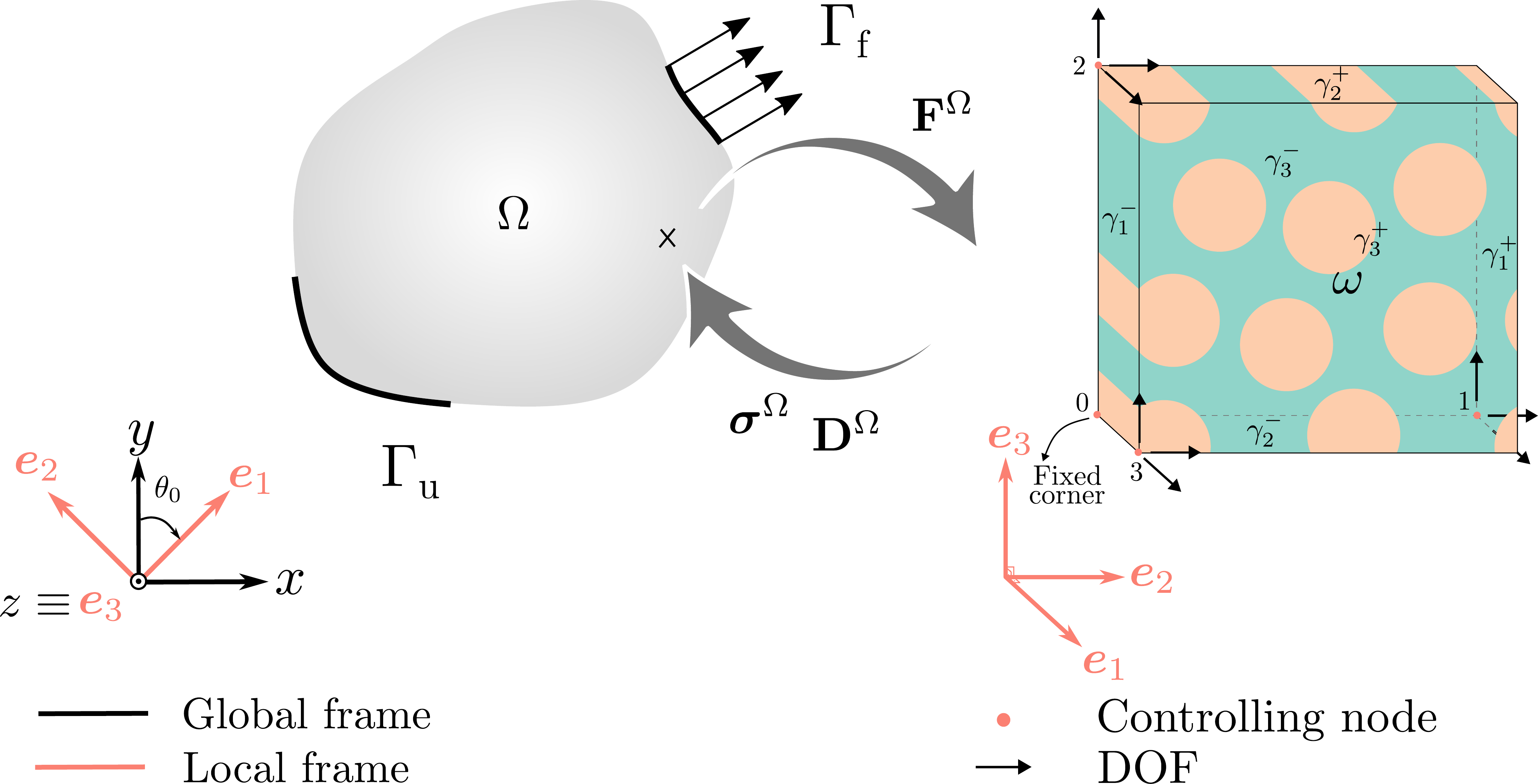}
\caption{Scheme of concurrent multiscale framework with two scales (macro $\Omega$ and micro $\omega$) for composite material with off-axis loading.}
\label{fig:multimicro}
\end{figure}

In this work, the updated Lagrangian formulation is adopted \cite{belytschko2014nonlinear}. The Cauchy stress and the deformation gradient are the stress and strain measures, respectively. At the microscale, regular constitutive models can be assigned to each of the phases. To solve the micromodel problem, periodic boundary conditions based on the macroscopic deformation gradient $\Fmacrotensor$ are usually employed. Further, if the local coordinate system of the RVE, in which fiber direction is always parallel to $\boldsymbol{e}_1$, is not aligned with the global one, we relate the two coordinate systems through a transformation matrix:
\begin{equation}
\mathbf{Q}_{0} \, (\theta_0) = \begin{bmatrix}
\cos (\theta_0) & -\sin(\theta_0) & 0\\
\sin(\theta_0) & \cos(\theta_0) & 0 \\
0 & 0 & 1
\end{bmatrix},
\label{eq:rotmatrix}
\end{equation}
where $\theta_0$ marks the angle between the $y$-axis and the fiber direction, as depicted in \Cref{fig:multimicro}. This matrix allows $\Fmacrotensor$ to be transformed from the global to the local coordinate system:
\begin{equation}
\label{eq:fromglobaltolocal}
\mathbf{F}_{\textrm{L }}^{\Omega} = \mathbf{Q}_{0} \Fmacrotensor \mathbf{Q}_{0}^{\textrm{T}}, 
\end{equation}
where the subscript ``L" refers to the local frame. Once the microscopic problem has converged, a computational homogenization procedure is performed to bridge the two scales:
\begin{equation}
\label{eq:macrostress}
\macrostressvec_{\textrm{L}} (\mathbf{X}^{\Omega}) = \frac{1}{|\omega|} \int_{\omega} \microstressvec (\mathbf{X}^{\omega}) \, \textrm{d}\omega.
\end{equation} 
Similarly, to transform the stress from the local to the global coordinate system, we use:   
\begin{equation}
\label{eq:fromlocaltoglobal}
\boldsymbol{\sigma}^{\Omega} = \mathbf{Q}_0^{\textrm{T}}\boldsymbol{\sigma}^{\Omega}_{\textrm{L}} \mathbf{Q}_0.
\end{equation}  

The formulation is completed with the definition of the macroscopic constitutive tangent $\mathbf{D}^{\Omega}$. For that, automatic differentiation, perturbation methods based on finite differences or condensation procedures can be used depending on the memory allocation and computational efficiency requirements \cite{TEMIZER2008tangent,Zhu2020}. In this work, our reference solutions are experiments and we do not perform the full-order FE$^2$ due to the exceedingly high computational cost. The main bottleneck, in this case, comes from the coupling between the two equilibrium problems, where the solution of the macroscopic displacement field defines the boundary condition for the RVEs, which in turn provide the missing homogenized constitutive model, requiring iterative solves of a large number of microscopic FE problems. In \Cref{sec:alternatives} we elaborate on the alternatives to this approach explored in this work.   

\subsection{Constitutive models}
\label{subsec:constmodels}
In this section, we discuss the two constitutive models used in the composite RVE adopted in the applications. These are exactly the same models as in \cite{Maiaetal2024} and the reference works \cite{Kovacevic2022,Kovacevic2022exp,Kovacevic2024creep}, and therefore we skip their derivation and only summarize their main features and the novelties explored in this work. The first constitutive model is a hyperelastic transversely isotropic model based on the formulation by \cite{Bonet1998} with slight modifications from \cite{Kovacevic2022}, which is assigned to the fibers. 

For the matrix, we use the Eindhoven Glassy Polymer (EGP), a rate and path-dependent elasto-viscoplastic model. 		The Cauchy stress from this model consists of three contributions. The first is a hydrostatic one, which is hyperelastic, the second is a hardening part that represents polymer chain reorientation and the last part comes from the driving stress, which is the component that allows for complex behavior to be incorporated and is rate and path-dependent. This component is given by the sum of several Maxwell elements, also known as modes, connected in parallel, and define a Non-newtonian flow for plastic deformations. In each mode, the spring is associated with a shear modulus and the dashpot to a viscosity function that depends on the stress, temperature and aging. It is also possible to consider multiple relaxation processes.

\subsection{Alternatives to the full multiscale formulation}
\label{sec:alternatives}
\begin{figure}[!h]
\centering
\includegraphics[width=\textwidth]{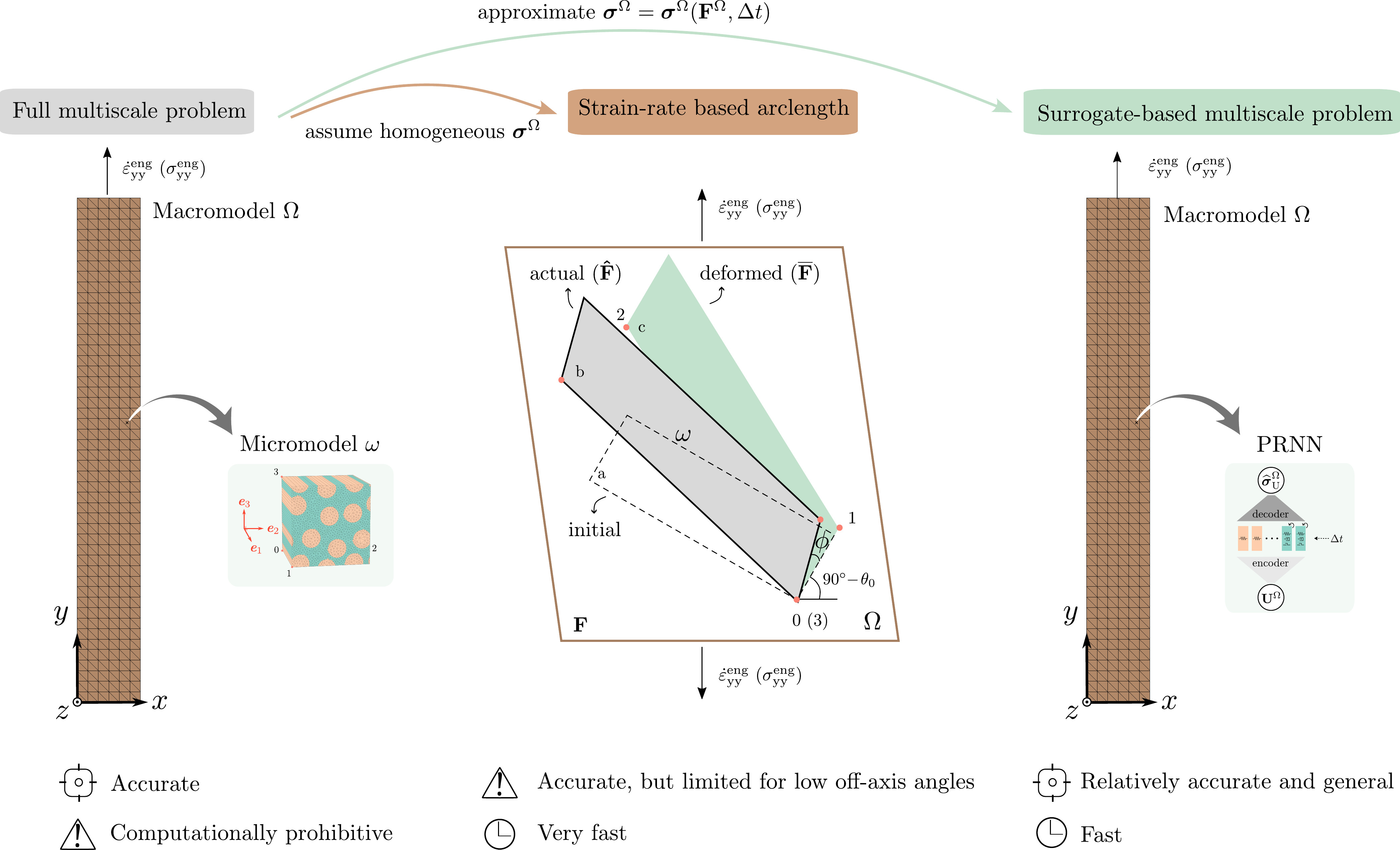}
\caption{Alternative approaches to replace full multiscale problem for modeling composites under off-axis loading.}
\label{fig:alternatives}
\end{figure}

In \Cref{fig:alternatives}, two alternatives to the full multiscale formulation are presented. Results obtained with both are compared in this paper. Before introducing the proposed surrogate-based multiscale approach, the single-scale micromechanical approach is summarized in this section. The single-scale model consists of a custom constant strain-rate arclength from \cite{Kovacevic2022exp}, which can be adapted to a force-controled framework to creep analysis \cite{Kovacevic2024creep}.  Although highly efficient and accurate in many scenarios, the assumptions made to simplify the entire specimen macroscopic domain to a single macroscopic point have limited validity. 

For results presented in this paper, slight changes are implemented with respect to the equations shown in the original formulation for the constant strain-rate case \cite{Kovacevic2022} due to the different choice of stress measure. In \cite{Kovacevic2022} and \cite{Kovacevic2022exp}, the Cauchy stress was used. Later, in a comparison with the creep experiments, the engineering stresses and strains were taken into account to agree with the engineering strain and stress reported in the experiments \cite{Kovacevic2024creep}. Here, we follow that choice and use the engineering quantities for comparison at all times. 

To apply a constant engineering strain-rate ($\dot{\varepsilon}_{\textrm{yy}}^{\textrm{eng}}$) under uniaxial engineering stress at the global frame ($\sigma_{\textrm{yy}}^{\textrm{eng}}$) on the macroscopic element, equivalent homogenized deformation and stress states need to be derived in the local frame. It is assumed that the specimen edge on which the loading is applied is prevented from rotating, while relative transverse translation is free, resulting in the parallelogram shape illustrated in the scheme depicted in \Cref{fig:alternatives}. A special constraint equation was derived by equating $F_\text{yy}$ to the value imposed from the input to ensure the global constant strain-rate condition while maintaining a uniaxial stress state in global frame.

Finally, to model the creep experiments, a similar framework based on the developments presented so far was proposed in \cite{Kovacevic2024creep}. The main changes regard the inclusion of the creep stress in the external force components for a force-controlled analysis and the consideration of the following condition:
\begin{equation}
\label{eq:creepstress}
\sigma_{\textrm{yy},\textrm{t}}^{\textrm{eng}} = \min \, \Bigl( \sigma_{\textrm{yy},t-1}^{\textrm{eng}} + \dot{\sigma}_{\textrm{yy}}^{\textrm{eng}} \Delta t, \sigma_{\textrm{yy}}^{\textrm{max}}\Bigr),
\end{equation} 
where $\dot{\sigma}_{\textrm{yy}}^{\textrm{eng}}$ is the engineering stress-rate at the loading phase before reaching the maximum stress level $\sigma_{\textrm{yy}}^{\textrm{max}}$. 

\subsection{Physically Recurrent Neural Network (PRNN)} 
\label{subsec:prnn}

In contrast to the microscale-based approach described previously, in this alternative we replace the RVE with a PRNN trained on stress-strain snapshots of it (see \Cref{fig:alternatives}). We summarize in \Cref{alg:prnnmacroroutine} the main coordinate system transformations needed to handle global and local frames when evaluating the homogenized stress response. To compute the macroscopic tangent stiffness matrix, we use the chain rule: 
\begin{equation}
\label{eq:tangentmacropoint}
\mathbf{D}^{\Omega} = \pd{\predstress}{\predstresslocal} \pd{\predstresslocal}{\Fmacrolocal} \ \pd{\Fmacrolocal}{\Fmacrotensor},
\end{equation} 
where the second term corresponds to the tangent obtained from the network: 
\begin{equation}
\label{eq:tangentprnn}
\pd{\predstresslocal}{\Fmacrolocal} = \pd{\predstresslocal}{\predstressunrot} \pd{\predstressunrot}{\Umacrolocal} \ \pd{\Umacrolocal}{\Fmacrolocal}  + \ \pd{\predstresslocal}{\Rmacrolocal} \ \pd{\Rmacrolocal}{\Fmacrolocal}. 
\end{equation}
 
\begin{algorithm}
    \SetKwInOut{Input}{Input}
    \SetKwInOut{Output}{Output}
    \Input{homogenized deformation gradient at global frame $\Fmacrotensor$, time increment $\Delta t$, initial off-axis angle $\theta_0$}
    \Output{homogenized stress $\widehat{\boldsymbol{\sigma}}^{\Omega}$, tangent stiffness matrix $\mathbf{D}^{\Omega}$}
    compute transformation matrix $\mathbf{Q}_0$ according to \Cref{eq:rotmatrix} based on $\theta_0$ \\
    transform strain from global to local frame: $\Fmacrotensor_{\textrm{L}} \leftarrow \mathbf{Q}_{0} \Fmacrotensor \mathbf{Q}_0^{\textrm{T}}$ \\
    perform polar decomposition on strain: $\mathbf{R}^{\Omega}_{\textrm{L}}, \mathbf{U}_{\textrm{L}}^{\Omega} \leftarrow$ \texttt{polarDecomposition}($\Fmacrotensor_{\textrm{L}}$) \\ 
    apply PRNN to stretch tensor: $\widehat{\boldsymbol{\sigma}}^{\Omega}_{\textrm{U}} \leftarrow $  \texttt{PRNN}($\Umacrotensor_{\textrm{L}}, \Delta t$) \\
    retrieve stress in the original coordinate system at the local frame: $\widehat{\boldsymbol{\sigma}}^{\Omega}_{\textrm{L}} \leftarrow \mathbf{R}_{\textrm{L}}^{\textrm{T}} \widehat{\boldsymbol{\sigma}}^{\Omega}_{\textrm{U}} \mathbf{R}_{\textrm{L}}^{\Omega}$ \\   
    transform stress from local to global frame: $\widehat{\boldsymbol{\sigma}}^{\Omega} \leftarrow \mathbf{Q}_{0}^{\textrm{T}} \widehat{\boldsymbol{\sigma}}^{\Omega}_{\textrm{L}} \mathbf{Q}_0$ \\
    compute tangent stiffness matrix $\mathbf{D}^{\Omega}$ according to \Cref{eq:tangentmacropoint} \\
   \textbf{return} ($\widehat{\boldsymbol{\sigma}}^{\Omega}$, $\mathbf{D}^{\Omega}$)
   \caption{Evaluation of macroscopic integration point using the PRNN}
   \label{alg:prnnmacroroutine}
\end{algorithm}

\begin{figure}[!h]
\centering
\hfill
\subfloat[Compact representation of PRNN\label{fig:sparsedecoder}]{\includegraphics[width=.35\textwidth]{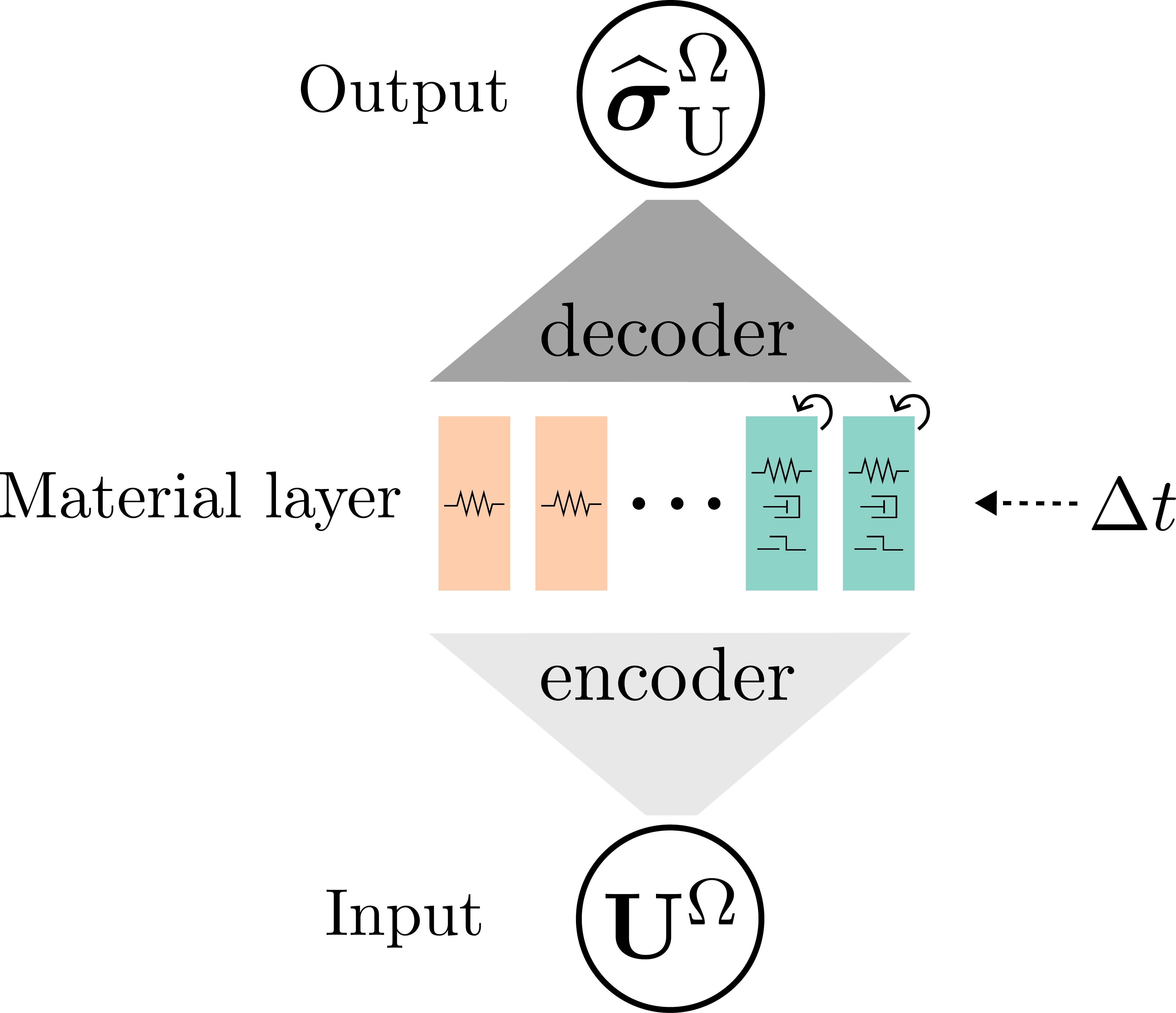}}
\hfill
\subfloat[Fictitious material point with history-dependent constitutive model unrolled in time
\label{fig:ficmatpoint}
]{\includegraphics[width=.49\textwidth]{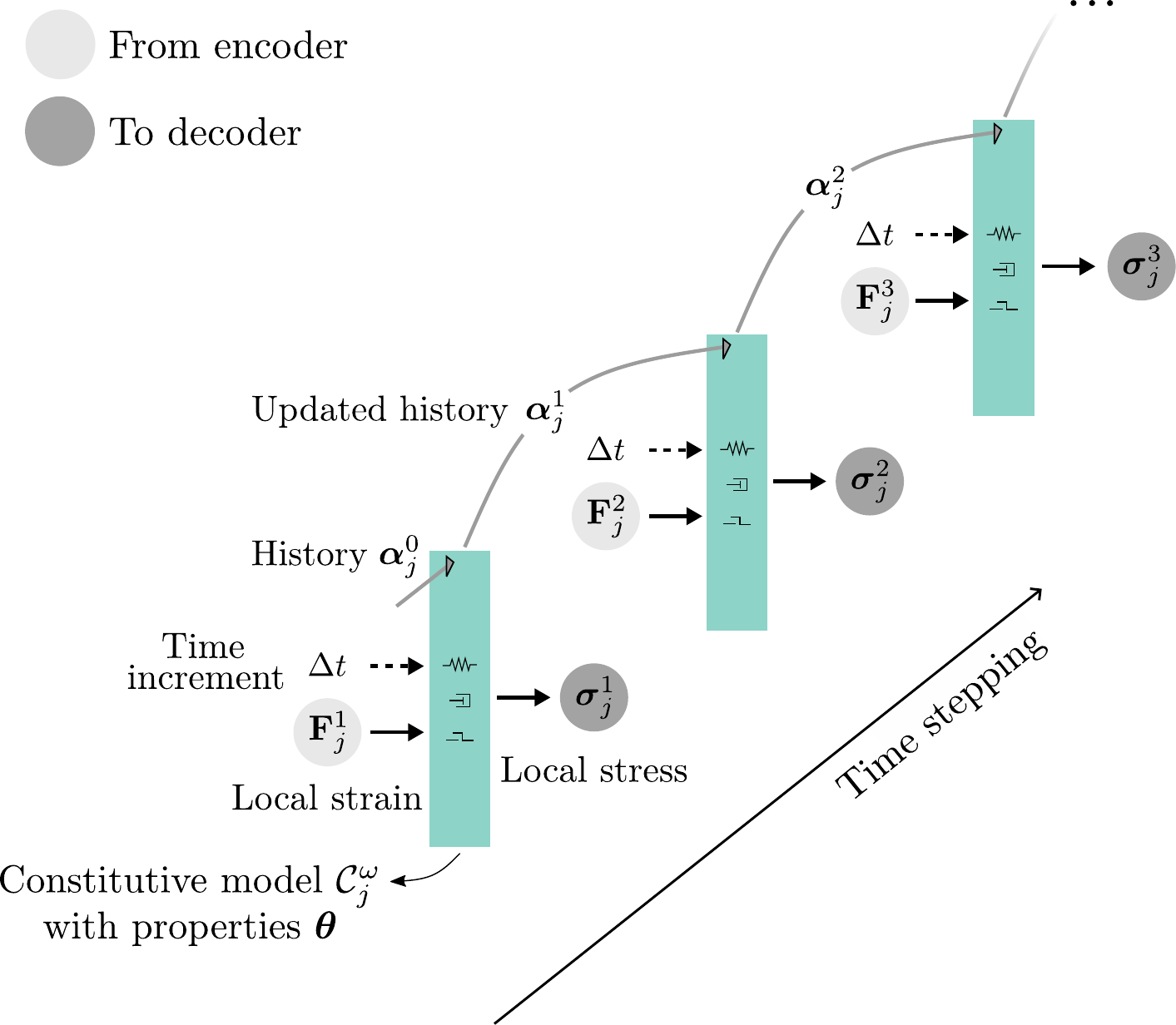}}
\hfill \
\caption{Architecture of PRNN and fictitious material point unrolled in time.}
\end{figure}

For details on implementation, training aspects or architectural choices, the reader is directed to \cite{Maiaetal2024}. Here, we highlight the main features of the network. In the PRNN, we preserve the constitutive models used in the micromodel, as well as their material properties, and embed them in an encoder-decoder architecture, as illustrated in \Cref{fig:sparsedecoder}. Through the encoder, we learn a set of values that we interpret as the local strain of \textit{fictitious microscopic material points}. Then, the constitutive model associated to each point is used to compute stresses and, in case of a history-dependent model, updated internal variables. Storing the internal variables from one time step to another, as illustrated in \Cref{fig:ficmatpoint} for a rate and path-dependent model $\mathcal{C}_j^{\omega}$, allows history dependence to arise naturally. Compared to the full-order model, the encoder learns a task analogous to the solution of the original microscale equilibrium problem, and the decoder corresponds to the stress homogenization. 

Having constitutive models in the network also raises a few constraints that regular NNs do not need to meet. For example, since the output of our encoder are deformation gradients, these must have positive determinant. For that purpose and taking into account the zero strain-stress state ($\Umacrotensor = \mathbf{I} \rightarrow \predstresslocal = \mathbf{0}$), we propose the weights from the encoder to be a symmetric matrix per material point. This arrangement avoids negative determinants and only requires 6 learnable parameters per material point. For the decoder, we employ a sparse architecture, in which only the matching stress components contribute to the macroscopic stress. Therefore, including encoder and decoder, each material point is associated with 12 learnable parameters. 

As for the constitutive model itself, both models present in the RVE, \material{matrix} and \material{fiber}, were embedded in the material layer. The splitting ratio between them is a hyper-parameter, but here we follow with the proportion from our previous contribution: the hyperelastic and elasto-viscoplastic models correspond to \qty{25}{\%} and \qty{75}{\%} of the material points, respectively, rounding the number of hyperelastic models up when the total number of points is even but not divisible by 4. Finally, in the following section, we present the transfer learning strategy used in this work to extrapolate in terms of material properties.
\begin{figure}[!h]
\centering
\includegraphics[width=\textwidth]{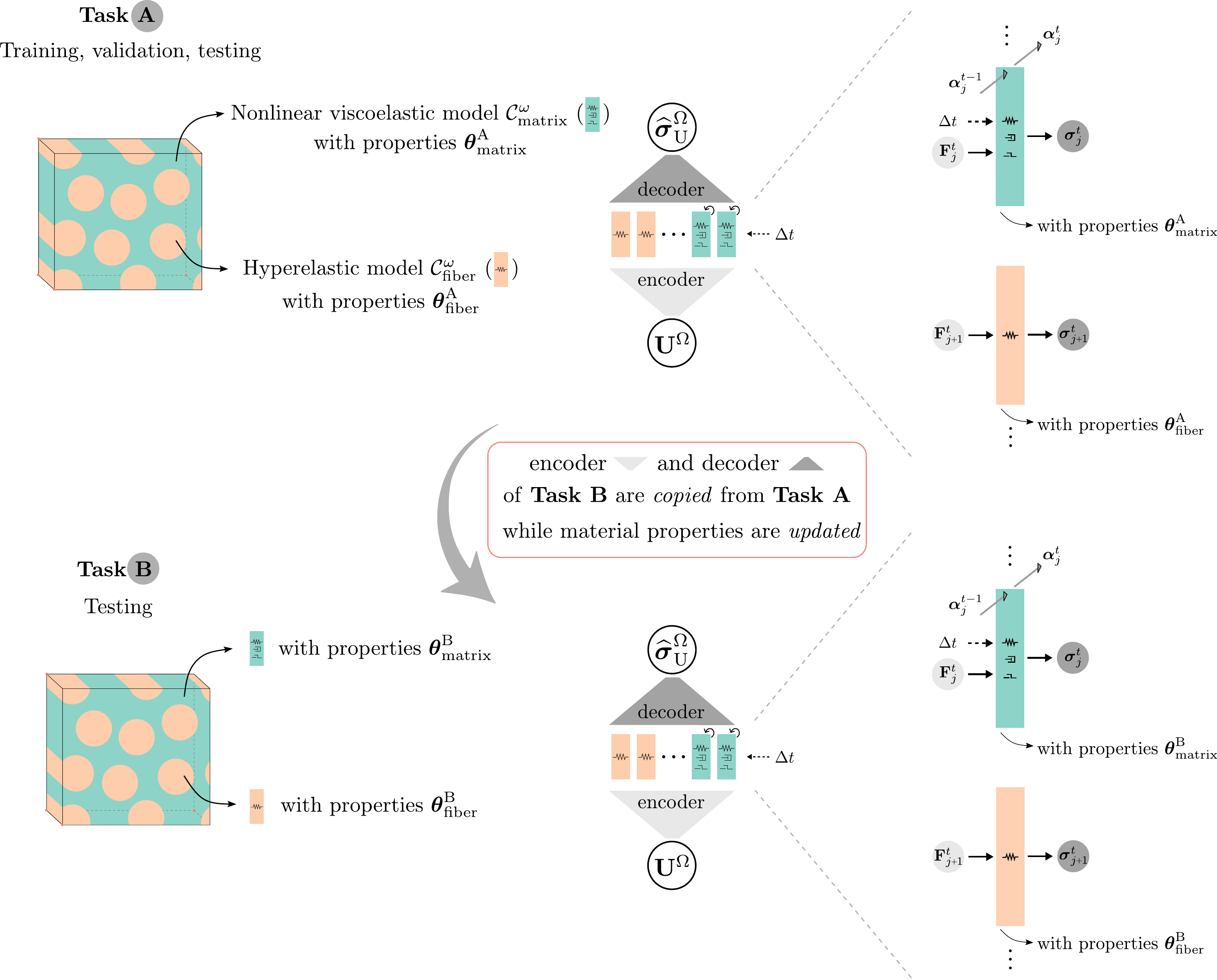}
\caption{From training with micromodel with a set of properties $\boldsymbol{\theta}_{(\cdot)}^{\textrm{A}}$ to extrapolating to micromodel with set of properties $\boldsymbol{\theta}_{(\cdot)}^{\textrm{B}}$ without retraining.}
\label{fig:tl}
\end{figure}

\subsection{Transfer learning}
\label{subsec:tl}
In this paper, we take a step further and explore how one can leverage the design choices in the material layer to transfer from one set of properties to another. The motivation for that comes from the fact that each mode in \material{matrix} is associated to 15 internal variables, to account for the full relaxation spectrum, a total of 240 internal variables would be needed per material point evaluated by this constitutive model. Such large sum of internal variables combined with the use of finite differences to compute the derivatives related to the backpropagation in time would slow down training significantly. 

The solution to that if using regular data-driven surrogate models would involve either training from scratch with the full 16 modes or use some transfer learning strategy, such as warm start. In both cases, additional computational effort would be necessary to (re)train the surrogates. Here, because material properties in the constitutive model are embedded in the material layer, we keep the network trained on a single mode and update the material properties in the network without any retraining. This scheme is illustrated in \Cref{fig:tl} with the models used in this work, but the approach is general and applicable to other constitutive models. The ability to extrapolate to different material properties in the online phase is yet another benefit of having an explainable function in our latent space. 

\section{Constant strain-rate experiments}
\label{sec:csrexp}

In this section, we update the material properties to transfer from a PRNN with single-mode \material{matrix} model to a multi-mode one. Then, we present the results obtained with the FEPRNN for modeling the constant strain-rate experiments. For this, a displacement controlled analysis was carried out, with the displacements in all three directions being constrained at the bottom and top surfaces of the coupon, except for the $y$ direction at the top, where the specimen is pulled at a constant strain-rate. We employ the adaptive stepping scheme from \cite{vanderMeer2012} to deal with potential convergence issues. In addition to that, we include a brief study on the use of oblique end-tabs as a possible adaptation to the experiments for a more uniform strain distribution. 

\subsection{From single to multi-mode PRNN}
\label{subsec:tlcsr}
To assess the performance of the networks from \cite{Maiaetal2024} trained with single-mode properties in extrapolating to the RVE with the full relaxation spectrum \setcsr{matrix}, we consider the largest training set from that work and keep the same optimal layer size, that is, 144 proportional GP-based curves and 8 material points in total (from which 6 are evaluated by \material{matrix} and 2 by \material{fiber}). In this particular application, the update in the material properties offers a further possibility for maximizing efficiency in the transfer learning scheme. Instead of considering the full relaxation spectrum at once, we investigate the gradual mode addition in the online phase. The idea comes from the fact that the latest modes have increasingly small contributions to the elastic regime, and could, in principle, be left out from the PRNN without loss of accuracy, leading to an even faster model evaluation. The gradual mode addition study is illustrated in \ref{sec:appendixmodes}, where we found that 5 modes are enough to accurately predict the response of RVEs generated with the full relaxation spectrum. Therefore, for the remainder of this section, we use 5 modes in all FEPRNN simulations.   

\subsection{Comparison with experimental results}
\begin{figure}[!ht]
\centering
\subfloat[Varying boundary conditions in micromodel for $\theta_0$ = \qty{15}{\degree}\label{fig:15degcsrwandwoutphi}]{
\includegraphics[width=0.5\textwidth]{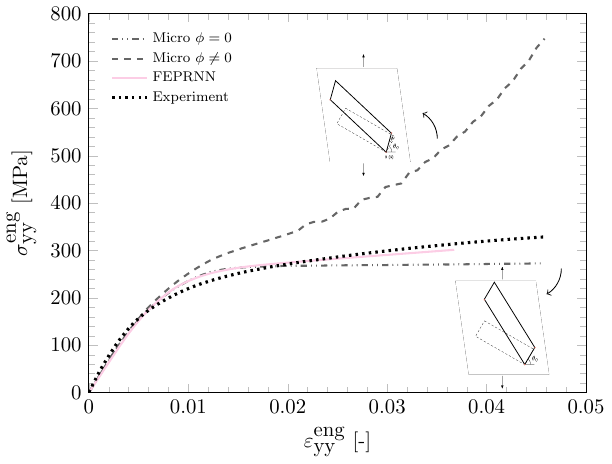}}
\subfloat[Different off-axis angles and micromodel boundary conditions \label{fig:vardegcsrwandwoutphi}]{\includegraphics[width=0.5\textwidth]{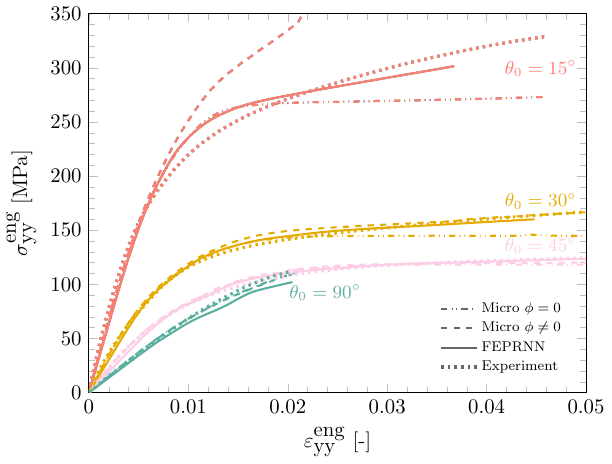}}
\caption{Stress-strain curve for $\dot{\varepsilon}_\textrm{yy}^{\textrm{eng}} =10^{-4} \, \unit{\per\second}$ using different methods and boundary conditions for multiple off-axis angles.}
\label{fig:difanglescsr}
\end{figure}

We start with the worst case scenario in \cite{Kovacevic2022exp}: off-axis angle $\theta_0 = \qty{15}{\degree}$ and strain-rate $\dot{\varepsilon}_\textrm{yy}^{\textrm{eng}} = 10^{-4} \, \unit{\per\second}$. In \Cref{fig:15degcsrwandwoutphi}, we plot the custom arclength response, referred as ``Micro", with and without the fiber rotation update $\phi$ against the experimental stress-strain curve. We also show results from the PRNN-based approach to the multiscale problem, named as ``FEPRNN". For this, we use the mesh shown in \Cref{fig:alternatives}, consisting of 576 (linear) wedge elements integrated with 1 point, and the PRNN discussed in \Cref{subsec:tlcsr} as the homogenized constitutive model. 

\begin{figure}[!h]
\centering
\includegraphics[width=.525\textwidth]{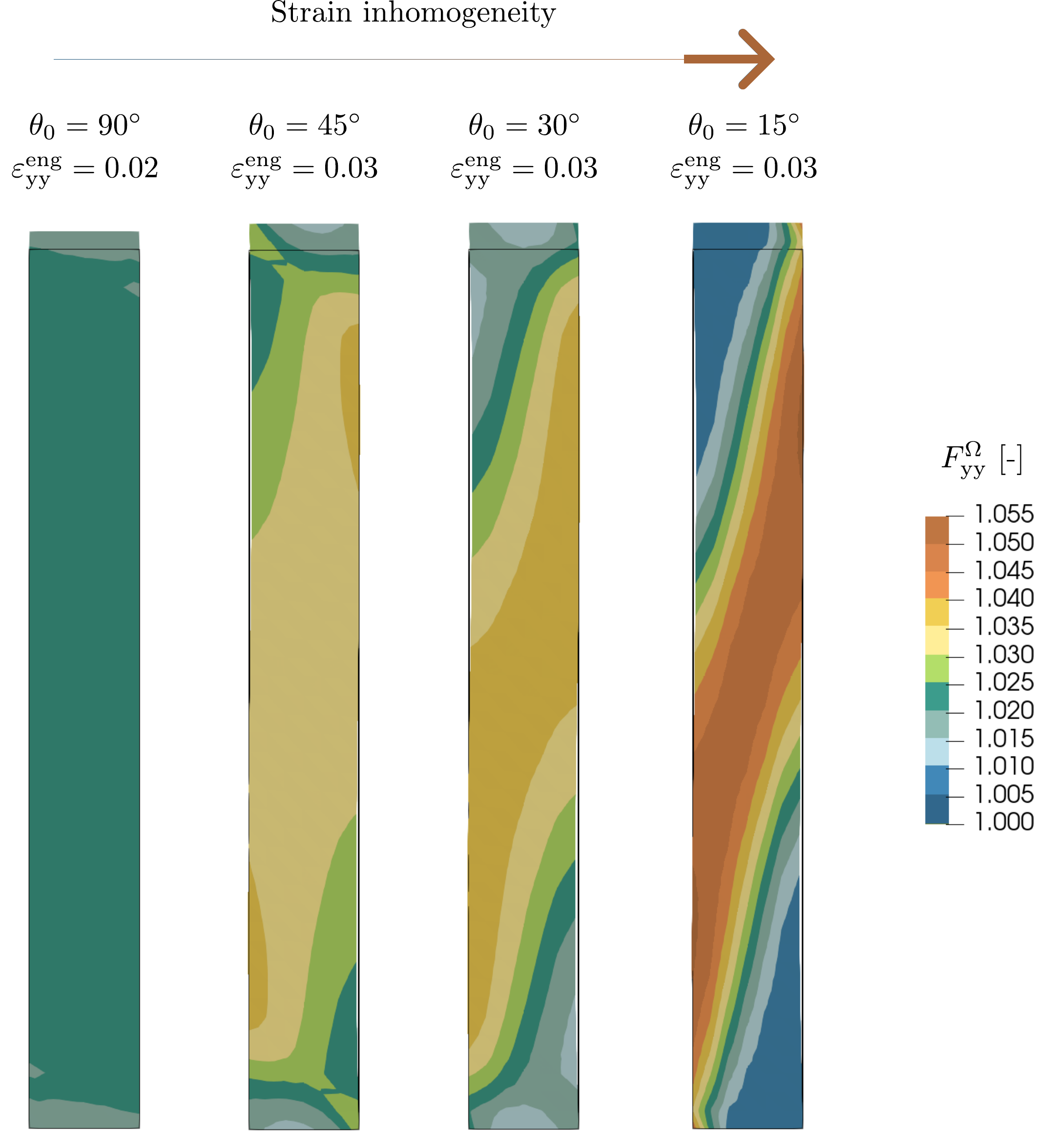}
\caption{Macroscopic strains from FEPRNN simulations for $\dot{\varepsilon}_{\textrm{yy}}^{\textrm{eng}} =10^{-4} \, \unit{\per\second}$ and different off-axis angles.}
\label{fig:fullfieldstraincsr}
\end{figure}

\Cref{fig:vardegcsrwandwoutphi} shows an overview of how the FEPRNN and micromodel-based solutions compare with the experiment for different off-axis angles. Comparing the responses using the micromodel with and without the reorientation angle $\phi$ in \Cref{fig:vardegcsrwandwoutphi}, it becomes clear how its effect depends heavily on the initial off-axis angle. Though these are not new results, it is only now possible to verify 1) how well the reorientation angle computed for the micromodel solution matches with the average angle computed from the PRNN-based simulation and 2) how well the assumption of macroscopic homogeneity holds. For this end, we show in \Cref{fig:fullfieldstraincsr} the full-field of strains in the loading direction for different off-axis angles. Note how the variation in the macroscopic strain becomes larger as the off-axis angle $\theta_0$ decreases. Clearly, the assumption behind the micromechanical analysis, that the specimen is in a macroscopically uniform state, is no longer valid. 

\begin{figure}[!ht]
\centering
\hfill
\begin{minipage}[b]{.41\linewidth}
\subfloat[$\theta_0 = \qty{45}{\degree}$]{\label{fig:phideg45}
\includegraphics[width=.95\textwidth]{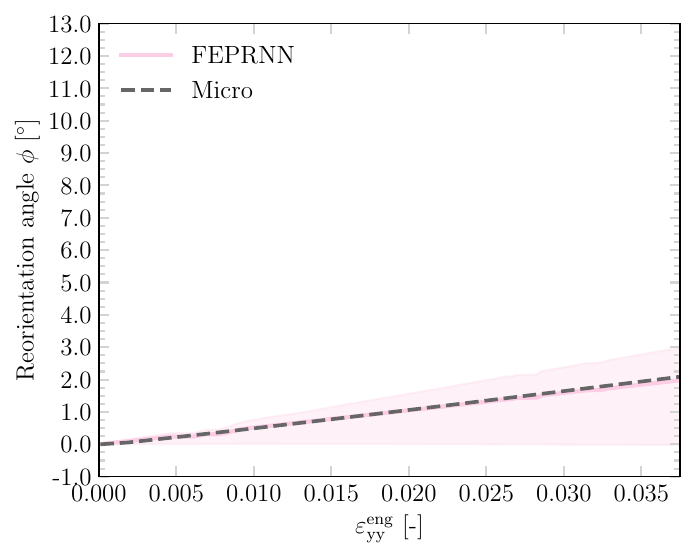}}
\hfill
\subfloat[$\theta_0 = \qty{15}{\degree}$]{\label{fig:phideg15}
\includegraphics[width=0.95\textwidth]{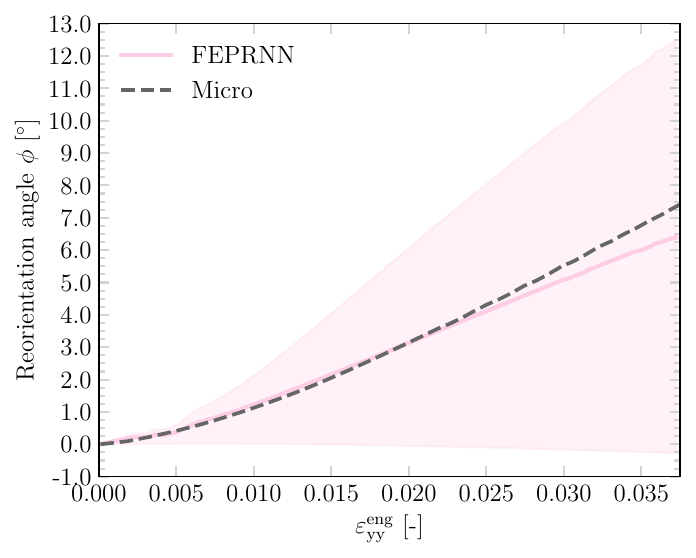}}
\end{minipage}
\hfill
\begin{minipage}[b]{.33\linewidth}
\subfloat[Snapshots of full field of angle $\phi$ using PRNN]{\label{fig:fullfieldphi}\includegraphics[width=1\textwidth]{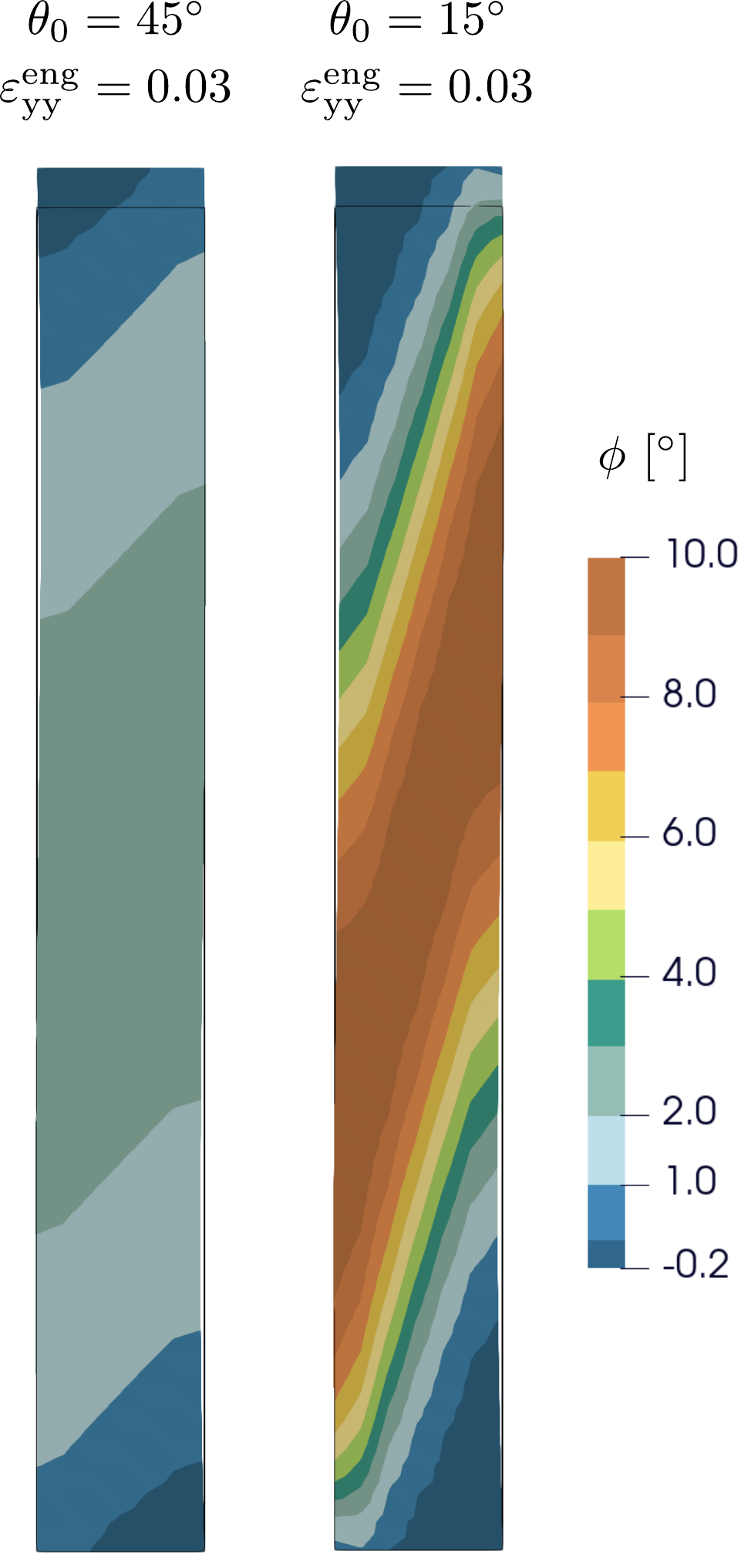}}
\end{minipage}
\hfill \ 
\caption{Reorientation angle for different initial off-axis angles for constant strain-rate experiment with $\dot{\varepsilon}_{\textrm{yy}}^{\textrm{eng}} =10^{-4} \, \unit{\per\second}$. Shaded areas in a) and b) represent the envelope of maximum and minimum values over the entire macroscopic domain, while solid lines correspond to the average.}
\label{fig:phi15and45}
\end{figure}

To further illustrate the problem, we calculate the reorientation angle $\phi$ at each of the macroscopic integration points with the surrogate-based multiscale approach. For this, we make a parallel with the deformation states in \Cref{fig:alternatives} and imagine that there is an equivalent shape of $\hat{\mathbf{F}}^{\Omega}$, whose side $0\--1$ is not moving in direction 2, and therefore $\overline{F}^{\Omega}_{21} = 0$. From that condition, the angle $\phi$ can be expressed as:
\begin{equation}
\label{eq:phimultiscale}
\phi = \arctan(\hat{F}^{\Omega}_{21}/\hat{F}^{\Omega}_{11})
\end{equation}
where $\hat{F}^{\Omega}_{21}$ and $\hat{F}^{\Omega}_{11}$ are components of the deformation gradient in the local frame (see \Cref{eq:fromglobaltolocal}). Based on this expression, we compute the simple average of $\phi$ over the entire macroscopic domain for the off-axis angles $\qty{45}{\degree}$ and $\qty{15}{\degree}$ and compare it with the one calculated by the custom arclength model, as shown in \Cref{fig:phi15and45}. In addition to the average, the envelope corresponding to the highest and lowest angles over the entire specimen are plotted (the pink shaded areas). Note how the angle variation is dramatically higher in the \qty{15}{\degree} case. For both angles, the highest values are located at the center of the specimen (see \Cref{fig:fullfieldphi}) and the lowest ones, close to zero, near the grips, where movement is restricted. 

In the $\qty{45}{\degree}$ case, the mean reorientation angle follows remarkably well the angle computed from the arclength model, which is also translated in the good visual agreement seen in \Cref{fig:difanglescsr}. For the $\qty{15}{\degree}$ composite, the envelope of the angle $\phi$ becomes increasingly larger and the single micromodel loses representativeness for the average response of the specimen, more strongly in the averaged stress-strain curve of \cref{fig:15degcsrwandwoutphi} than in the mean orientation angle in \cref{fig:phideg15}. 

\begin{figure}[!h]
\centering
\includegraphics[width=.55\textwidth]{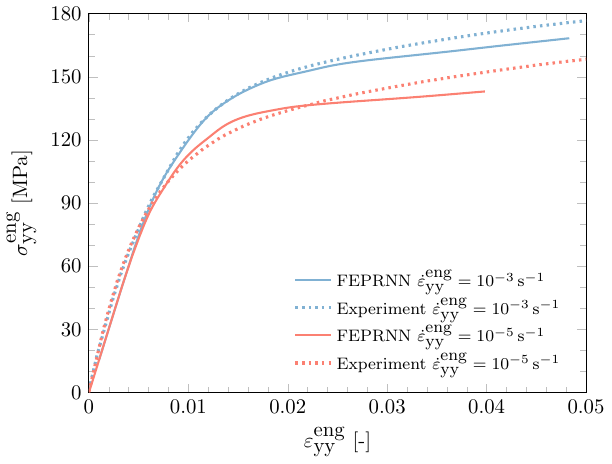}
\caption{Stress-strain curves for off-axis angle $\theta_0 = \qty{30}{\degree}$ and different strain-rates $\dot{\varepsilon}_\textrm{yy}^{\textrm{eng}}$.}
\label{fig:30degdifsrcsr}
\end{figure}

In addition to that, the multiscale simulations showed the presence of shear stress, suggesting that the boundary conditions in the experiments are not only far from the assumptions in the arclength model, but also from its initial goal, a uniaxial tensile test. Finally, we illustrate the generality of the framework with two other strain rates for the $\qty{30}{\degree}$ case: $10^{-3} \, \unit{\per\second}$ and $10^{-5} \, \unit{\per\second}$. \Cref{fig:30degdifsrcsr} shows the good match of the two engineering stress-strain curves obtained using the PRNN as the constitutive model with respect to the experimental curves. In both cases, the maximum error is around \qty{5}{\percent}, at the last time step. 

\subsection{Investigation into alternative setups}
\begin{figure}[!h]
\centering
\includegraphics[width=.09\textwidth]{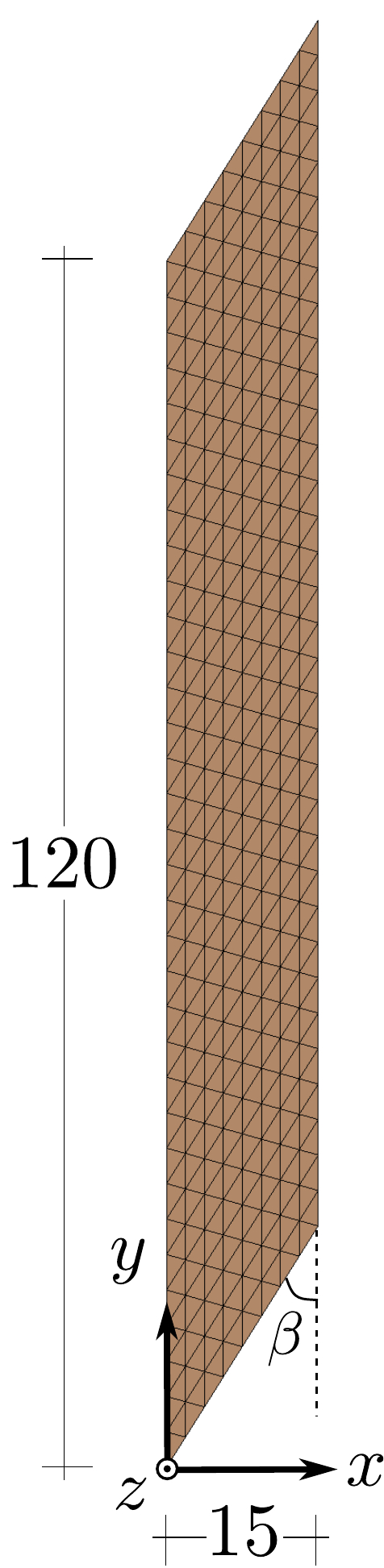}
\caption{Alternative coupon geometry with oblique-end tabs.}
\label{fig:obliqueangle}
\end{figure}

With the multiscale framework, we can also investigate the influence of the specimen geometry and boundary conditions on the response, which could be interesting in case one seeks to calibrate material models based on engineering stresses and strains coming directly from test frame readings. First, we analyze an alternative geometry with oblique end-tabs. The tab design proposed by \citet{SUN1993} is illustrated in \Cref{fig:obliqueangle}, with $\beta$ given by:
\begin{equation}
\beta = \text{acot} \, (-S_{16}/S_{11})
\label{eq:obliqueangle}
\end{equation}
where $S_{16}$ and $S_{11}$ are the compliance coefficients with respect to the global coordinate system. For linear elasticity, this design ensures a macroscopically homogeneous stress states under uniaxial loading. In the nonlinear case, the compliance coefficients change from one time step to another, but we adopt the values obtained at the undeformed configuration with 5 modes. 

\begin{figure}[!ht]
\centering
\hfill
\subfloat[\label{fig:obliqueaxial}Axial stress]{
\includegraphics[width=0.48\textwidth]{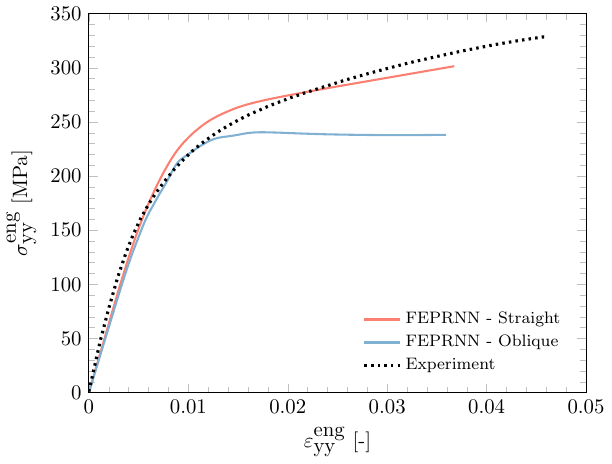}}
\hfill 
\subfloat[\label{fig:obliqueshear}Shear stress]{
\includegraphics[width=0.48\textwidth]{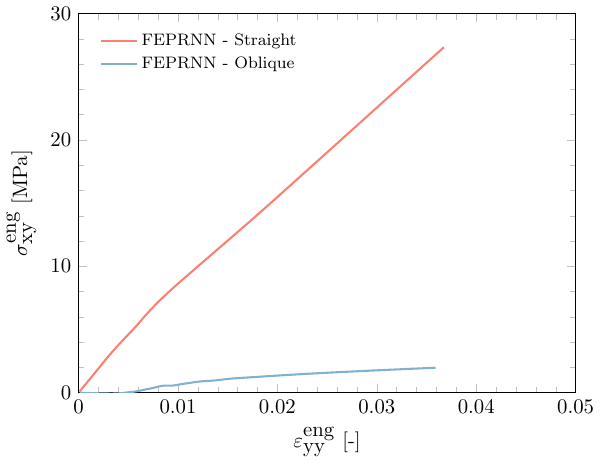}}
\hfill \
\caption{Stress-strain curves for $\theta_0 = \qty{15}{\degree}$ and $\dot{\varepsilon}_\textrm{yy}^{\textrm{eng}} =10^{-4} \, \unit{\per\second}$ using straight vs oblique end-tabs.}
\label{fig:obliquestraightengcurves}
\end{figure}

To illustrate the effectiveness of the oblique end-tabs in creating a nearly uniform stress-state, we continue with the most challenging case explored so far: off-axis angle of $\qty{15}{\degree}$ and strain-rate of $10^{-4} \, \unit{\per\second}$. For comparison, we plot the stress-strain curves for $\sigma_{\textrm{yy}}^{\textrm{eng}}$ and $\sigma_{\textrm{xy}}^{\textrm{eng}}$ the two types of end-tabs in \Cref{fig:obliquestraightengcurves}. Note that the oblique design is effective in reducing the shear almost to zero. The new geometry also results in a lower engineering stress level in the loading direction (axial).     

\begin{figure}[!h]
\centering
\includegraphics[width=.6\textwidth]{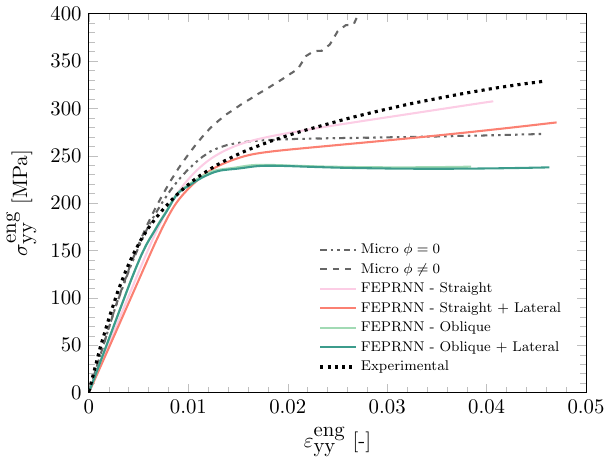}
\caption{Stress-strain curves for different end-tab designs and boundary conditions for $\theta_0 = 15^{\circ}$ and $\dot{\varepsilon}_\textrm{yy}^{\textrm{eng}} =10^{-4} \, \unit{\per\second}$.}
\label{fig:summarycsr}
\end{figure}

Another way of eliminating shear stress is to allow lateral movement of the specimen. \Cref{fig:summarycsr} shows the comparison between the two types of end-tabs and different boundary conditions with respect to the experiments. For the oblique design, allowing the lateral movement has little to no impact on the outcome since the shear stress is already quite limited in that case. On the other hand, with the straight end-tabs, allowing the lateral movement reduces the engineering stress by about \qty{10}{\percent}.

\begin{figure}[!h]
\centering
\includegraphics[width=.575\textwidth]{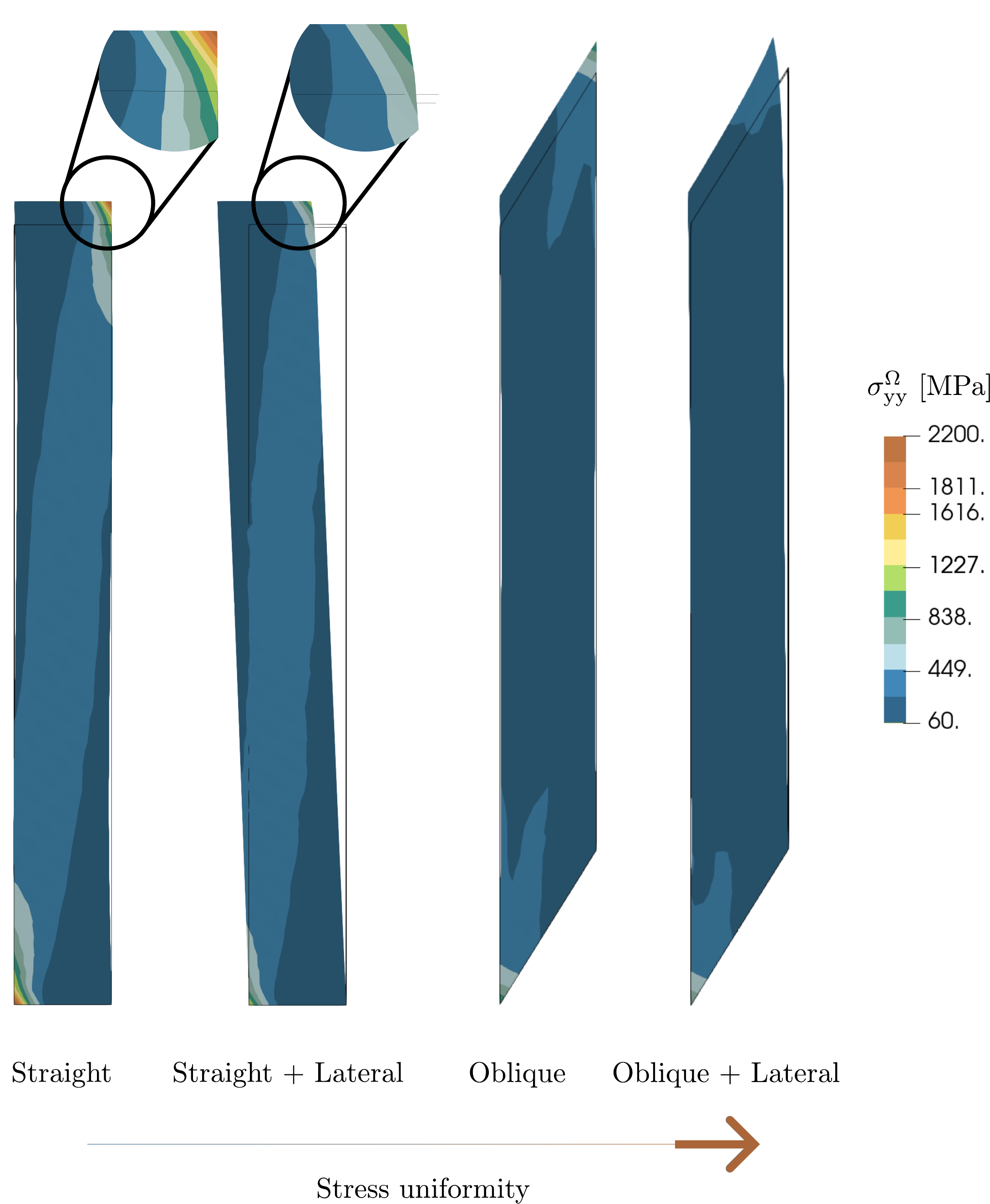}
\caption{Full field of stress predicted by FEPRNN at the same strain level ($\varepsilon_\textrm{yy}^{\textrm{eng}} = 0.03$) on specimen with straight and oblique end-tabs with and without lateral movement allowed for $\theta_0 = \qty{15}{\degree}$ and $\dot{\varepsilon}_\textrm{yy}^{\textrm{eng}} =10^{-4}\, \unit{\per\second}$.}
\label{fig:fullfieldstresscsr}
\end{figure}

Finally, the full-field stress distribution for the combinations discussed previously is shown in \Cref{fig:fullfieldstresscsr}. Note the stress concentration near the edges of the specimen and the higher stress band formed in the diagonal of the coupon, aligned with the fiber direction. These are still present in the (computational) experiment with lateral movement, but are significantly reduced with the use of oblique end-tabs on the right side.  

\section{Creep experiments}
\label{sec:creep}

In this section, we follow the same structure as in the previous section, only now for the creep experiments. Again, we employ the adaptive stepping scheme \cite{vanderMeer2012}, this time with larger bounds for the time increments as the creep simulations are significantly longer. We keep the same boundary conditions as in \Cref{sec:csrexp} now with a force controlled analysis according to \Cref{eq:creepstress}. 

\subsection{From one relaxation spectrum to another}
In \cite{Kovacevic2024creep}, a new relaxation spectrum for  \material{matrix} with a total of 2 processes and 8 modes was calibrated directly on the creep experiments. In these experiments, no extensometer was used, therefore measurements implicitly take into account the compliance effect of the machine grips. In addition to that, the shear modulus $G_{12}$ was changed from 45 MPa to 5 MPa for \material{fiber} \cite{Kovacevic2024creep}. Apart form these properties, since the same micromodel was used, the straightforward option would be to update the material properties in both \material{matrix} and \material{fiber} and test the PRNN performance with gradually increasing modes until convergence, in the same fashion as in \Cref{subsec:tlcsr}. This time, however, direct transfer to the new set of properties was not successful. To better understand the problem, we generated several test sets considering different combinations of matrix and fibers properties and found that the direct update of the material properties works, but on a certain range around the original properties used for training. This study can be found in \ref{sec:appendixdifga}. 

Hence, for modeling the creep experiments, we take a step back and train the PRNNs from scratch with the new set of properties \setcreep{fiber} and the first 2 modes (the first from each process) from \setcreep{matrix}. For the training, the same approach and optimization parameters as in \cite{Maiaetal2024} were used with two modifications: 1) when creating the proportional GP-based curves, the time increment in each curve is sampled from a log-uniform distribution $\Delta t\sim U \, (10^{-3} \, \unit{\s}, 10^3 \, \unit{\s})$, and 2) the maximum number of epochs was increased to 2000. Then, we apply the same approach as discussed in \Cref{subsec:tl} and \Cref{subsec:tlcsr} to go from a network trained on 2 modes from \setcreep{matrix} to one that predicts the response of a micromodel with the full relaxation spectrum. As this follows the same story-line as \Cref{subsec:tlcsr}, we leave the model selection and the gradual mode addition discussion for the \ref{sec:trcreep}. In short, we select the best network with 6 material points (4 evaluated by \material{matrix} and 2 by \material{fiber}) trained on 72 curves and tested on 150 curves of the same loading type and material properties. For the remaining FEPRNN simulations, we consider 6 modes in the fictitious points evaluated by \material{matrix}.

\subsection{Comparison with experimental results}
\label{subsec:resultscreep} 
\begin{figure}[!h]
\centering
\subfloat[$\theta_0 = \qty{15}{\degree}$, $\sigma_\textrm{yy}^\textrm{max}$ = 290 MPa\label{fig:creep15deg}]{\includegraphics[width=.49\textwidth]{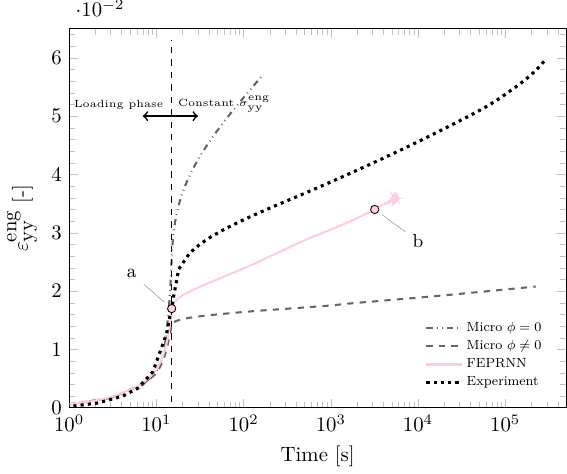}}
\hfill
\subfloat[$\theta_0 = \qty{90}{\degree}$, $\sigma_\textrm{yy}^\textrm{max}$ = 97 MPa\label{fig:creep90deg}]{\includegraphics[width=.505\textwidth]{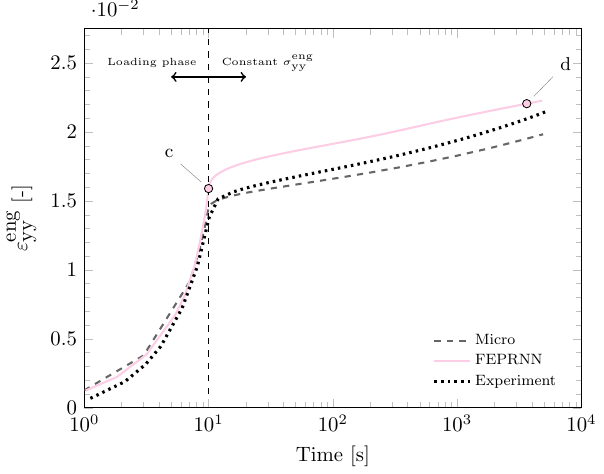}}
\hfill
\subfloat[Macroscopic strain field predicted by FEPRNN\label{fig:fullfieldcreep15and90}]{\includegraphics[width=.52\textwidth]{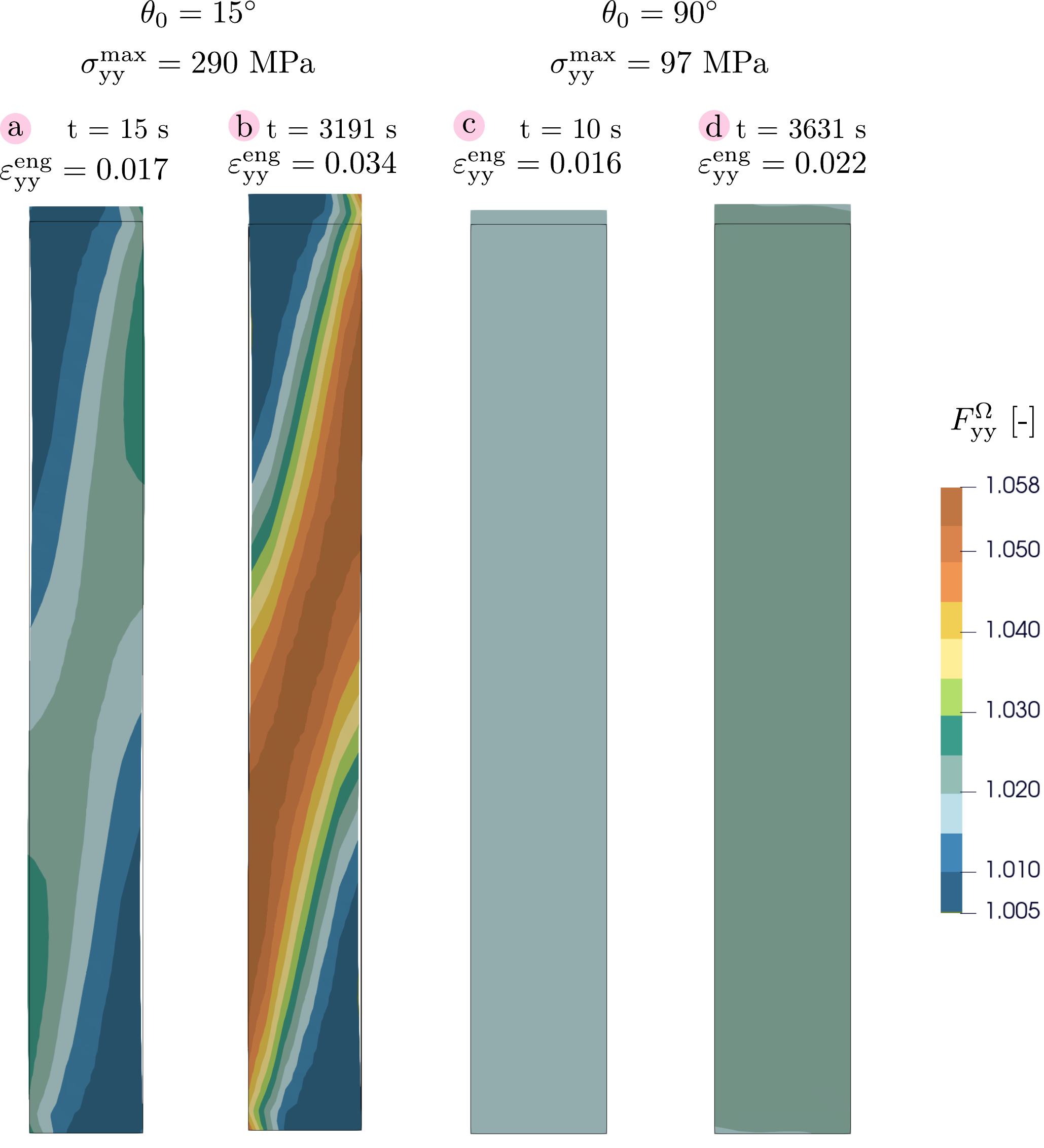}}
\caption{Strain-time curves and macroscopic strain distributions for off-axis angles \qty{15}{\degree} and \qty{90}{\degree}.}
\label{fig:creep15and90deg}
\end{figure}

We start the discussion with the $\theta_0 = \qty{15}{\degree}$ and $\theta_0 = \qty{90}{\degree}$ cases. \Cref{fig:creep15and90deg} shows the strain-time FEPRNN response for the two combinations of off-axis angle and engineering stress, as well as the micromodel results from \cite{Kovacevic2024creep}, without the consideration of micro-cracking. Similar to the constant strain-rate experiments, for $\theta_0 = \qty{15}{\degree}$, the effect of the reorientation angle update $\phi$ on the micromodel results is pronounced and a highly inhomogenenous strain field is obtained with the FEPRNN approach, as illustrated in \Cref{fig:fullfieldcreep15and90}. Despite the offset with respect to the experimental results, the creep response of the FEPRNN again still lies between the two bounds of the micromodel solution, showing that the (surrogate-based) multiscale simulation is vital to account for the macroscopic variations that the micromodel-based approach cannot.  

\begin{figure}[!h]
\centering
\includegraphics[width=.57\textwidth]{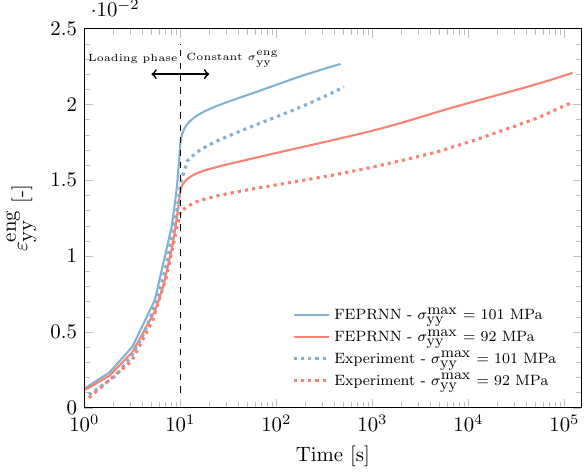}
\caption{Strain-time curves for $\theta_0 = \qty{90}{\degree}$ and different stress levels $\sigma_\textrm{yy}^\textrm{max}$.}
\label{fig:creepdifstress}
\end{figure}

Next, we investigate the scenario for which fiber reorientation is not relevant: the $\theta_0 = \qty{90}{\degree}$ case. This time, a uniform macroscopic strain distribution is obtained, as shown in \Cref{fig:fullfieldcreep15and90} with the snapshots of the macroscopic strain distribution corresponding to the points ``c" and ``d" marked in \Cref{fig:creep90deg}. In this case, we obtain satisfactory agreement of the FEPRNN with respect to the creep experiments. To further demonstrate the generality of the PRNN, we consider two different engineering stress levels. In \Cref{fig:creepdifstress}, a satisfactory agreement is achieved for the two stress levels with fixed $\theta_0 = \qty{90}{\degree}$, with maximum errors below \qty{10}{\percent}. 

\begin{figure}[!h]
\centering
\subfloat[$\theta_0 = \qty{30}{\degree}$, $\sigma_\textrm{yy}^\textrm{max}$ = 154 MPa\label{fig:creep30deg}]{
\includegraphics[width=.49\textwidth]{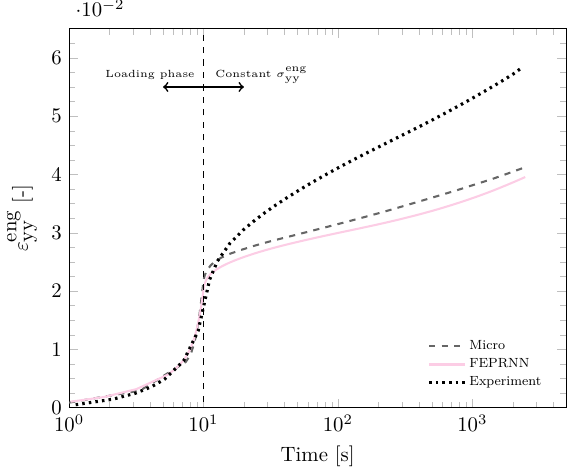}}
\hfill
\subfloat[$\theta_0 = \qty{45}{\degree}$, $\sigma_\textrm{yy}^\textrm{max}$ = 110 MPa\label{fig:45degdifbc}]{
\includegraphics[width=0.49\textwidth]{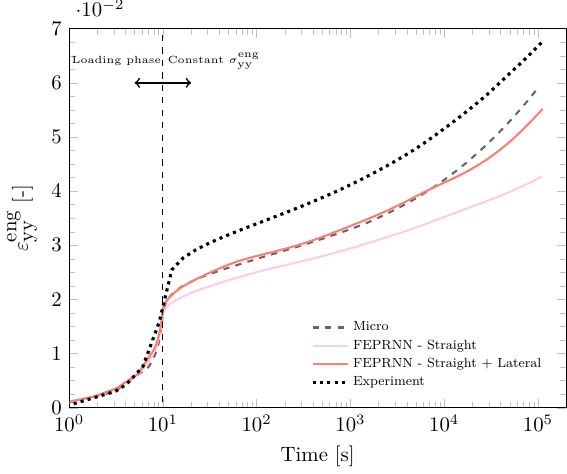}}
\hfill \
\caption{Strain-time curves and full-field snapshots for $\theta_0 = \qty{30}{\degree}$ and $\theta_0 = \qty{45}{\degree}$ and loading time t = $\qty{10}{\second}$.}
\label{fig:creepdifangle}
\end{figure}

For $\theta_0 = \qty{45}{\degree}$ and $\theta_0 = \qty{30}{\degree}$, the FEPRNN significantly undershoots the engineering strain, reaching at best the same response as the micromodel, as shown in \Cref{fig:creepdifangle}. In \Cref{fig:45degdifbc} specifically, the setup with more uniform macroscopic strain/stress distribution (straight + lateral) is the closest to the micromodel response and to the the experiments. This is the opposite of what is found in \Cref{sec:csrexp}, where the higher-fidelity reproduction of the experimental setup (straight) allowed better agreement with the experiments. One explanation lies in the fact that in that case, the material properties were calibrated on stress-strain curves obtained from extensometer measurements, while in the creep experiments no extensometer was used. As stated in \cite{Kovacevic2024creep}, the incorporation of more accurate results would demand the re-calibration of the relaxation spectrum used in \material{matrix}. With material properties that implicitly consider the machine grip effect, the multiscale simulation is correcting the same effect twice. 

\section{Conclusion}
\label{sec:conclusion}
In this paper, we address the challenges observed in \cite{Kovacevic2022exp, Kovacevic2024creep} to model constant strain-rate and creep experiments on unidirectional composites under off-axis loading. In those works, a simplified and efficient approach was proposed based on a single macroscopic point with special boundary conditions to reproduce an (assumed) uniaxial stress state. While successful in many cases, the framework showed limitations with lower off-axis angles. The main hypothesis is that the assumption of macroscopic uniform deformation did not hold in the experiments. To test this hypothesis with the same high-fidelity micromodel, a full multiscale approach would be necessary, which remains computationally prohibitive. However, with the development of the PRNN, a surrogate model recently tested on the same micromodel and constitutive models as in this work \cite{Maiaetal2024}, a surrogate-based multiscale analysis becomes within reach. Building on these contributions, it is now possible to verify the hypotheses raised in \cite{Kovacevic2022exp,Kovacevic2024creep} and to gain insights into the potential macroscopic strain/stress distribution in the experiments.

The PRNNs have a hybrid architecture that combines data-driven components and physics-based constitutive models to leverage better generalization properties with limited datasets. They have been used to model the homogenized response of micromodels with a range of constitutive laws, from linear elastic to elasto-viscoplastic models, including the ones considered in this work. With embedded models in the latent space, extrapolation is possible not only in terms of loading types but also in terms of material properties. Here, we harness this feature to simplify and reduce the computational cost of the offline phase by training with a fraction of the set of properties used in the references \cite{Kovacevic2022exp,Kovacevic2024creep}. Specifically, we generate and train the network with the first mode of each process in the elasto-viscoplastic model. Because each mode is associated with many internal variables, and these are at the core of the backpropagation in time, training with fewer modes translates into a more efficient process. We then show that these networks can be successfully transferred to our main task without any re-training effort by directly updating the material properties of the constitutive model. This feature is illustrated in the multiscale applications with the gradual mode addition. 

From the constant strain-rate experiments, and assuming that the good match between them and the multiscale simulations implies a faithful characterization of the former, we summarize the following findings:   
\begin{itemize}
\item Due to the tension/shear coupling, the constrained movement of the specimen by the machine grips, and the straight end-tabs, significant shear is present when the coupon is loaded in tension, except for the extreme cases $\theta_0 = \qty{0}{\degree}$ and $\theta_0 = \qty{90}{\degree}$;
\item Allowing free shear deformation of the coupon is not enough to achieve uniform strain and stress distributions. Although it does eliminate the shear stress measured at the top of the coupon, the stress concentration problem is still present (only slightly less so than in the fully constrained case);
\item We observe highly inhomogeneous macroscopic strain and stress fields for the lower off-axis angles, particularly for $\theta_0 = \qty{15}{\degree}$, with strain and stress concentrations at the corners near the grips of the specimen. In such cases, the large macroscopic variability indicates that the experimental results should not be  interpreted as direct uniaxial stress-strain measurements of the composite material;
\item We verify the hypotheses raised in the reference work of the micromechanical model \cite{Kovacevic2022exp} to explain the lack of fitting for $\theta_0 = \qty{15}{\degree}$. The experimental response is indeed a combination of material points that rotate to align with the loading direction (at the center of the specimen) and material points that have little to no rotation (near the surface of the end-tabs). For the remaining studied off-axis loading cases, there is moderate ($\theta_0 = \qty{30}{\degree}$ and $\theta_0 = \qty{45}{\degree}$) to negligible ($\theta_0 = \qty{90}{\degree}$) strain variation;
\item To achieve a more uniform stress-strain uniaxial state in these experiments, we investigate the use of oblique end-tabs. It is confirmed that these also help reduce the stress concentration near the grips and greatly reduce the manifestation of shear stress across the specimen. 
\end{itemize}

For the creep experiments, a new set of properties was adopted. This time, they were calibrated directly on the stress-strain experimental creep measurements obtained without the use of extensometers, which implicitly accounts for the machine grip effect. With that, we observe limited accuracy in the multiscale simulations for the intermediary off-axis angles ($\theta_0 = \qty{30}{\degree}$ and $\theta_0 = \qty{45}{\degree}$) compared to the micromodel solution. Our results highlight the pitfalls of calibrating the material parameters from experiments without reliable strain measurements and under the assumption of homogeneous stress fields

For the two remaining cases, opposite scenarios are observed, one where the strain distribution is quite uniform with no machine grip effect ($\theta_0 = \qty{90}{\degree}$) and the other with large macroscopic variability ($\theta_0 = \qty{15}{\degree}$). For the former, we obtain a relatively good fit at different engineering stress levels. Whereas for $\theta_0 = \qty{15}{\degree}$, despite the convergence issues and gap between experimental and FEPRNN results, the approximate response is in the middle of the two extreme cases of the micromodel-based approach and follows the same trend as the experiment. These findings, along with the insights from the constant-strain experiments, showcase the potential of PRNNs to enable efficient and robust multiscale analysis under variable loading conditions. 

\section*{Acknowledgments}
The authors acknowledge the TU Delft AI Initiative for their support through the SLIMM AI Lab. FM acknowledges financial support from the Dutch Research Council (NWO) under Vidi grant 16464.

\section*{Data availability}
A centralized repository for PRNN developments can be found at https:/github.com/\-SLIMM-Lab/pyprnn

\appendix
\section{Mode addition with properties calibrated for constant-strain rate experiments}
\label{sec:appendixmodes}
\begin{figure}[!ht]
\centering
\hfill
\subfloat[Uniaxial stress with $\Delta \lambda = \num{1e-4}$ and $\Delta t = \qty{1}{\second}$]{
\includegraphics[width=0.475\textwidth]{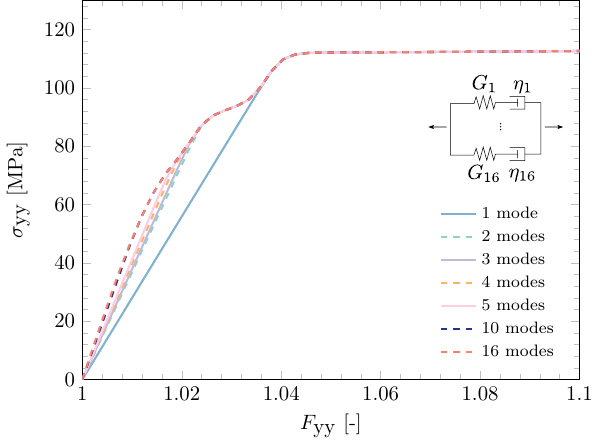}}
\hfill
\subfloat[Accuracy of PRNN trained on 144 proportional GP-based curves and a single mode over test set $\mathcal{T}_\textrm{16}$ considering different number of modes in the fictitious material points evaluated by \material{matrix}\label{fig:errorMultiModetl}]{
\includegraphics[width=0.5\textwidth]{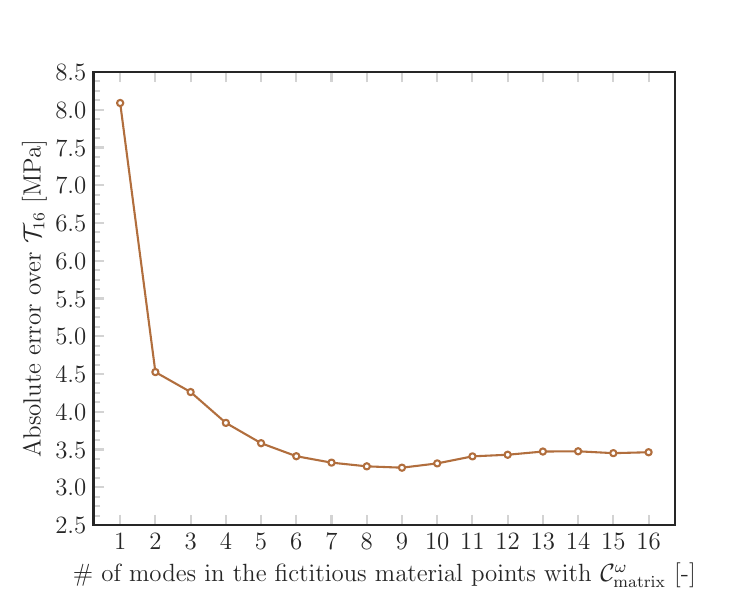}}
\hfill \
\caption{Stress-strain curves using \material{matrix} and PRNN accuracy with different number of modes from \setcsr{matrix}.\label{fig:egpmodes}}
\end{figure}

In \Cref{fig:egpmodes}, we show several stress-strain curves obtained using the EGP model (\material{matrix}) with different number of modes calibrated for constant-strain rate experiments \setcsr{matrix}. Although the pre-peak regime response can change significantly with more modes, the post-peak response (plastic regime) is the same. For testing the networks trained on single-mode from \cite{Maiaetal2024} in extrapolating to the response of a RVE with the full relaxation spectrum, we select the best network over test set $\mathcal{T}_{\textrm{1}}$, which consists of 150 curves with the same loading type and material properties as the ones used for training. This loading type is deemed to be representative of unseen loading conditions as it includes different unloading/reloading cycles. Then, on a new test set $\testunlgpbased$ with 150 curves generated using the 16 modes in \material{matrix}, we assess the accuracy of the network by updating the material properties considering a gradual addition of modes in the fictitious material points evaluated by \material{matrix}, as shown in \Cref{fig:errorMultiModetl}. The error reaches a minimum around 9 modes with \qty{5.9}{\percent}.

\section{Limitations on the transfer learning for creep experiments}
\label{sec:appendixdifga}
\begin{figure}[!ht]
\centering
\includegraphics[width=0.8\textwidth]{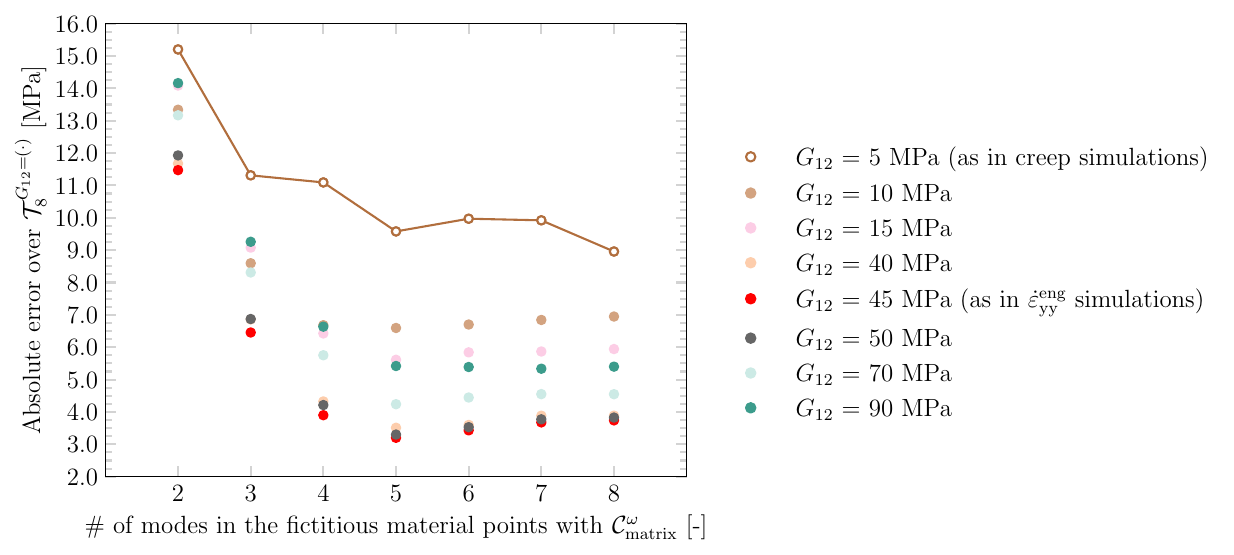}
\caption{Accuracy with gradual mode addition in the fictitious material points evaluated by \material{matrix} of PRNN trained on micromodel simulations with a single mode from $\boldsymbol{\theta}_{\textrm{matrix}}^{\dot{\varepsilon}_\textrm{yy}}$ over test sets with micromodel simulations with 8 modes from \setcreep{matrix} and different fiber shear moduli.}
\label{fig:errorMultiModeCreeptl}
\end{figure}

The following combination of material properties was considered to generate test sets of 150 proportional GP-based curves each: full relaxation spectrum calibrated for the creep experiments, referred as \setcreep{matrix}, and varying shear modulus $G_{12}$ in the material properties considered for the fibers \setcreep{fiber}. The only difference between \setcsr{fiber} and \setcreep{fiber} is the $G_{12}$, which changed from 45 MPa to 5 MPa. With the network parameters obtained training on the first mode of \setcsr{matrix} and \setcsr{fiber}, we update the material properties in the fictitious material points with \material{fiber} to match the corresponding shear modulus being tested, while the modes used in \material{matrix} were gradually updated with increasing number of modes from the new spectrum. Though the test sets have slightly different stress-strain ranges, the comparison in terms of absolute errors in \Cref{fig:errorMultiModeCreeptl} is meaningful and reveals a few insights. Firstly, going from the 1st mode of one relaxation spectrum (\setcsr{matrix}) to a multi-mode prediction of another (\setcreep{matrix}) is possible. Second, the best performance is achieved around 5 modes. And lastly, in this case, the change in the shear modulus is the most important parameter in determining how well we can extrapolate. For the value considered in the creep simulations ($G_{\textrm{12}} = 5$ MPa), the lowest error observed is still quite high, around 9 MPa. 

\section{Training and transfer learning for creep experiments}
\label{sec:trcreep}
\begin{figure}[!h]
\centering
\subfloat[Envelope of highest and lowest errors from 10 initialization of PRNNs trained on different number of proportional GP-based curves and material layer size over validation set $\mathcal{V}_{\textrm{2}}$\label{fig:creepmodelselection}]{\includegraphics[width=.5\textwidth]{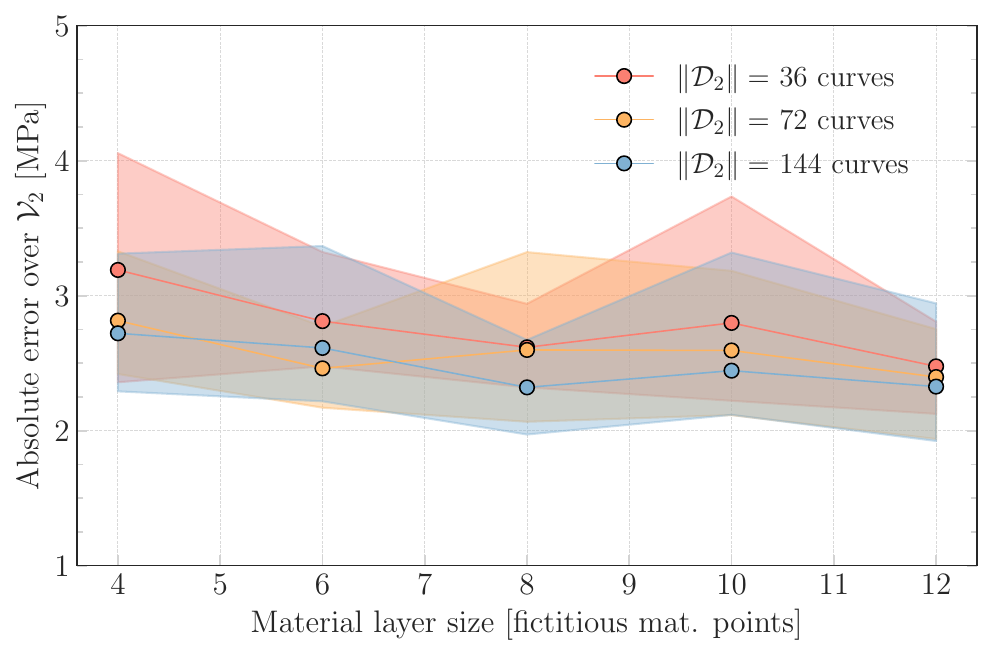}}
\hfill
\subfloat[Accuracy with mode addition in PRNN with 6 material points trained on the first two modes over test sets with the full relaxation spectrum $\mathcal{T}_{\textrm{8}}$ \label{fig:creepextrapolation}]{\includegraphics[width=.45\textwidth]{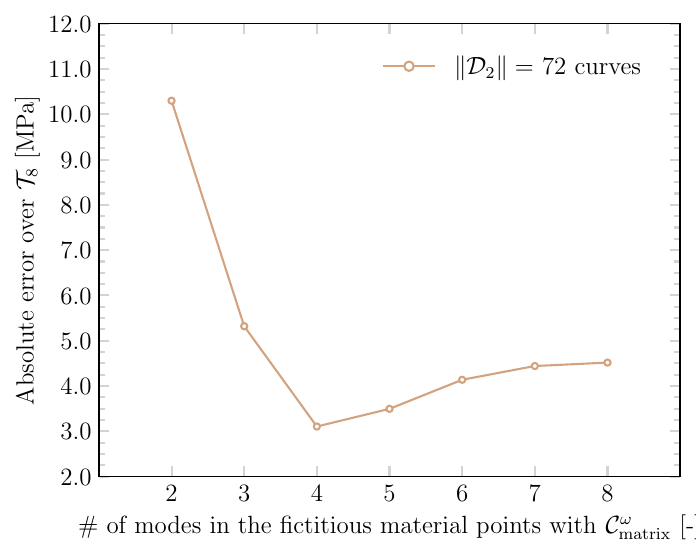}}
\caption{Model selection procedure for PRNNs trained on the micromodel with material properties calibrated on creep experiments.}
\label{fig:errorMultiModeCreepTrainedon2mo}
\end{figure}

\Cref{fig:creepmodelselection} shows the envelope of highest and lowest errors over a validation set $\mathcal{V}_2$ with 150 proportional GP-based curves for the different initialization for a number of material layer and training set sizes, from which we select the networks with 6 fictitious material points (4 evaluated by \material{matrix} and 2 by \material{fiber}) trained on 72 proportional GP-based curves. To assess the performance of the PRNNs in extrapolating to more modes, first we select the network with the lowest error over a test set $\mathcal{T}_{\textrm{2}}$ with 150 proportional GP-based curves and same material properties as in the training. Then, on a new test set $\mathcal{T}_{\textrm{2}}$ that considers RVEs with the new full relaxation spectrum, we test the accuracy of that network with a varying number of modes (see \Cref{fig:creepextrapolation}). In this case, the errors reach a minimum around \qty{3}{\MPa} ($\approx$ \qty{10}{\percent}) when 4 modes are considered. 

\biboptions{sort&compress}
\bibliographystyle{elsarticle-num-names} 
\bibliography{references.bib}

\end{document}